\documentclass{amsart}

\usepackage[all]{xy}

\usepackage{latexsym}

\usepackage{amssymb}
\usepackage{amsfonts}
\usepackage{amscd}
\usepackage{amsmath,amsthm}

\newtheorem{lemma}{Lemma}[section]
\newtheorem{proposition}[lemma]{Proposition}
\newtheorem{theorem}[lemma]{Theorem}
\newtheorem{corollary}[lemma]{Corollary}

\newtheorem{prop}[lemma]{Proposition}

\newtheorem{cor}[lemma]{Corollary}

{
}

\theoremstyle{definition}

\newtheorem{definition}[lemma]{\sl Definition}

\theoremstyle{remark}

\newcommand{\Hom}{\operatorname{Hom}}
\newcommand{\sHom}{{\mathcal{H}}{\it{om}}_{{\sf Gr }A}}
\newcommand{\usHom}{{\underline{{\mathcal{H}}{\it{om}}}}_{{\sf Gr }A}}
\newcommand{\sExt}{{\mathcal{E}xt}_{{\sf Gr }A}}
\newcommand{\usExt}{{\underline{\mathcal{E}xt}}_{{\sf Gr }A}}
\newcommand{\dlim}{\underset{n \rightarrow \infty}{\lim}}

\newcommand{\el}{\mathcal{L}}

\numberwithin{equation}{section}

\begin{document}

\title{An abstract characterization of noncommutative projective lines}
\author{A. Nyman}
\address{Department of Mathematics, 516 High St, Western Washington University, Bellingham, WA 98225-9063}
\email{adam.nyman@wwu.edu}
\keywords{}
\thanks{2010 {\it Mathematics Subject Classification.} Primary 14A22; Secondary 16S38}

\begin{abstract}
Let $k$ be a field.  We describe necessary and sufficient conditions for a $k$-linear abelian category to be a noncommutative projective line, i.e. a noncommutative $\mathbb{P}^{1}$-bundle over a pair of division rings over $k$.  As an application, we prove that $\mathbb{P}^{1}_{n}$, Piontkovski's $n$th noncommutative projective line, is the noncommutative projectivization of an $n$-dimensional vector space.
\end{abstract}

\maketitle

\tableofcontents

\pagenumbering{arabic}

\section{Introduction}
Throughout this paper, we work over a field $k$.

\subsection{Motivation}   In noncommutative algebraic geometry, the primary objects of study are abelian categories which bear similarities to categories of coherent sheaves over schemes.  When $k$ is algebraically closed, a classification of categories like the category of coherent sheaves over a smooth projective curve has been achieved \cite{artinstafford}, \cite{vdbreiten}.  In the case that $k$ is not algebraically closed, there is a general theory of noncommutative curves due to Kussin \cite{kussin2}, and those of genus zero have been studied extensively \cite{kussin}.

In this paper, we study a notion of noncommutative projective line which specializes the concept of noncommutative $\mathbb{P}^{1}$-bundle due to M. Van den Bergh \cite{vandenbergh}.  To motivate it, recall that the ordinary projective line over $k$ can be constructed as the projectivization of a two-dimensional vector space $V$, i.e. $\mathbb{P}^{1} = \mathbb{P}(V) = \operatorname{Proj} \mathbb{S}(V)$, where $\mathbb{S}(V)$ is the symmetric algebra of $V$.  In order to construct a noncommutative analogue of $\mathbb{P}^{1}$, we replace $V$ by a (suitably well-behaved) bimodule $M$ over a pair of division rings, we replace $\mathbb{S}(V)$ by Van den Bergh's noncommutative symmetric algebra $\mathbb{S}^{nc}(M)$ (see Section \ref{section.ncsym}), and we replace the category of coherent sheaves over $\mathbb{P}(V)$ by the appropriate quotient category of graded coherent right $\mathbb{S}^{nc}(M)$-modules (see below), which we denote by $\mathbb{P}^{nc}(M)$.  Although one has to make some restrictions on $M$ so that the category $\mathbb{P}^{nc}(M)$ is homologically well behaved (see the statement of Theorem \ref{theorem.supermain}(5)), once these restrictions are made, $\mathbb{P}^{nc}(M)$ behaves in many respects like the commutative projective line.  If, for example, $M$ has left-right dimension $(2,2)$ and we let $\mathbb{P}^{nc}(M)$ denote the category of finitely generated graded right $\mathbb{S}^{nc}(M)$-modules modulo those that are eventually zero, then the category $\mathbb{P}^{nc}(M)$ is noetherian and has many properties in common with the commutative projective line \cite{newchan}.  Even if $M$ does not have left-right dimension $(2,2)$, the space $\mathbb{P}^{nc}(M)$ can still be defined but is usually not noetherian.  However, the results of this paper and \cite{chancoherent} establish that $\mathbb{P}^{nc}(M)$ is hereditary, satisfies Serre duality, and has the property that each of its objects is a direct sum of its torsion part and a sum of line bundles.  Motivated by these results, we refer to spaces of the form $\mathbb{P}^{nc}(M)$ as {\it noncommutative projective lines}.  Our main result (Theorem \ref{theorem.supermain}(4) and (5)) is a description of necessary and sufficient conditions on a $k$-linear abelian category to be a noncommutative projective line.

Noncommutative projective lines are related to other notions of noncommutative curve.  For example, recall that Kussin defines a noncommutative curve of genus zero to be a small $k$-linear noetherian, abelian, Ext-finite category $\sf C$ with a Serre functor inducing the relevant form of Serre duality, an object of infinite length, and a tilting object.  The curve ${\sf C}$ is {\it homogenous} if for all simple objects $\mathcal{S}$ in ${\sf C}$, $\operatorname{Ext}_{{\sf C}}^{1}(\mathcal{S},\mathcal{S}) \neq 0$.  Concrete examples of such curves are the noncommutative projective lines $\mathbb{P}^{nc}(M)$ where $D_{0}$ and $D_{1}$ are division rings which are finite dimensional over $k$, and $M$ is a $k$-central $D_{0}-D_{1}$-bimodule of left-right dimension $(1,4)$ or $(2,2)$.  In \cite[Theorem 3.10]{tsen} and \cite[Theorem 3.1]{witt}, the following converse is established:

\begin{theorem} \label{theorem.tsen}
Every homogeneous noncommutative curve of genus zero has the form $\mathbb{P}^{nc}(M)$ for some division rings $D_{0}$ and $D_{1}$ and some $k$-central $D_{0}-D_{1}$-bimodule $M$ of left-right dimension $(1,4)$ or $(2,2)$.
\end{theorem}

The notion of noncommutative projective line we work with is broad enough to encompass examples not described using Kussin's definition.  In particular, both generic fibers of noncommutative ruled surfaces not finite over their centers (\cite{patrick}, \cite{newchan}), and Piontkovski's noncommutative projective lines (\cite{piont}), are noncommutative projective lines in our sense.
The first class of examples is relevant to Artin's conjecture \cite{problems}, which loosely states that noncommutative surfaces which are infinite over their center are birationally ruled.  In this context, our main theorem specializes to give conditions under which a category is equivalent to the generic fiber of a noncommutative ruled surface.  The second class of examples implies that many noncommutative projective lines are non-noetherian, and leads us to work more generally within the framework of coherent noncommutative algebraic geometry \cite{polish}.  In this framework, a noncommutative projecive line is a quotient category of the form
$$
\mathbb{P}^{nc}(M) := {\sf coh }\mathbb{S}^{nc}(M)/{\sf tors }\mathbb{S}^{nc}(M)
$$
where ${\sf coh }\mathbb{S}^{nc}(M)$ denotes the full subcategory of graded right $\mathbb{S}^{nc}(M)$-modules consisting of coherent modules and ${\sf tors }\mathbb{S}^{nc}(M)$ is the full subcategory of ${\sf coh }\mathbb{S}^{nc}(M)$ consisting of right-bounded modules.

\subsection{The main theorem in the finite over $k$ case}
We now describe a special case of our main result, the so-called finite over $k$ case, deferring a complete statement to Section \ref{section.props} (Theorem \ref{theorem.supermain}).  Its statement involves a collection of properties which may or may not hold for an arbitrary sequence
$$
\underline{\mathcal{L}} := (\mathcal{L}_{i})_{i \in \mathbb{Z}}
$$
of objects in an abelian $k$-linear category ${\sf C}$.  The properties we are interested in are as follows:  for all $i \in \mathbb{Z}$,
\begin{itemize}
\item{} $\operatorname{End }\mathcal{L}_{i}= :D_{i}$ is a division ring finite-dimensional over $k$, and $\operatorname{dim }_{k}D_{i}=\operatorname{dim }_{k}D_{i+2}$,

\item{} $\Hom(\el_{i},\el_{i-1})=0=\operatorname{Ext}^{1}(\mathcal{L}_{i},\mathcal{L}_{i-1})$,

\item{} Let $l_{i}=\operatorname{dim}_{D_{i+1}}\operatorname{Hom }(\mathcal{L}_{i},\mathcal{L}_{i+1})$.  Then there is a short exact sequence
$$
0 \rightarrow \mathcal{L}_{i} \rightarrow \mathcal{L}_{i+1}^{ l_{i}} \rightarrow \mathcal{L}_{i+2} \rightarrow 0
$$
such that the $l_{i}$ morphisms defining the left arrow are linearly independent over $D_{i+1}$.

\item{} $\operatorname{Ext}^{1}(\mathcal{L}_{i}, \mathcal{L}_{j})=0$ for all $j \geq i$.

\item{} for all $\mathcal{M}$ in ${\sf C}$, $\operatorname{Hom}(\mathcal{L}_{i}, \mathcal{M})$ is a finite-dimensional $D_{i}$-module.

\item{} $\underline{\mathcal{L}}$ is ample (see Section \ref{section.ampleness} for the definition of ampleness).
\end{itemize}
For example, if ${\sf C}$ is the category of coherent sheaves over $\mathbb{P}^{1}$, and $\mathcal{O}$ denotes the structure sheaf on $\mathbb{P}^{1}$, then the sequence $(\mathcal{O}(i))_{i \in \mathbb{Z}}$ satisfies these conditions.

Finally, we need to introduce the following terminology regarding bimodules.  We say a bimodule over a pair of division rings {\it is of type $(m,n)$} with $m,n$ nonnegative integers if it has left-dimension $m$ and right-dimension $n$ or left-dimension $n$ and right-dimension $m$.

Our main theorem specializes to the following:

\begin{theorem} \label{theorem.ft}
(Finite over $k$ case) The category ${\sf C}$ has a sequence $\underline{\mathcal{L}}$ satisfying the above six conditions if and only if
$$
{\sf C }\equiv \mathbb{P}^{nc}(M)
$$
where $M$ is a $D_{0}-D_{1}$-bimodule of type $(m,n) \neq (1,1), (1,2), (1,3)$ such that $\mathbb{S}^{nc}(M)$ is graded right-coherent.
\end{theorem}
The fact that noncommutative symmetric algebras are graded right-coherent has been proven under very general hypotheses \cite{chancoherent}.

We call Theorem \ref{theorem.ft} the finite over $k$ case as the rings $\operatorname{End }\mathcal{L}_{i}$ are assumed to be finite-dimensional over $k$.  This hypothesis does not hold for generic fibers of noncommutative ruled surfaces, and is not necessary in general.  The cost of removing the hypothesis is that we must then add other hypotheses to equate the left-dimension of a bimodule over division rings with the right-dimension of its right-dual.

As a consequence of Theorem \ref{theorem.ft} we show that the noncommutative projective lines, $\mathbb{P}^{1}_{n}$, $n \geq 2$, studied in \cite{piont} and \cite{kgroup}, are noncommutative projective lines in our sense.  To describe this application in more detail we recall the definition of $\mathbb{P}^{1}_{n}$.  To this end, we need the following result of J.J. Zhang \cite[Theorem 0.1]{zhang}: every connected graded regular algebra of dimension two over $k$ generated in degree one is isomorphic to an algebra of the form \begin{equation} \label{eqn.algebraa}
A=k \langle x_{1}, \ldots, x_{n} \rangle/ (b)
\end{equation}
where $b=\sum_{i=1}^{n}x_{i}\sigma(x_{n-i})$ for some graded automorphism $\sigma$ of the free algebra.  Piontkovski shows that such rings are graded right-coherent \cite[Theorem 4.3]{piont}, and defines $\mathbb{P}^{1}_{n}$ to be ${\sf cohproj }A$, where $A$ is as above, and ${\sf cohproj }A$ denotes the category of coherent graded right $A$-modules modulo those that are eventually zero.  Piontkovski then proves \cite[Theorem 1.5]{piont} that this category depends only on $k$ and $n$, and that such categories share many homological properties with the category of coherent sheaves on (commutative) $\mathbb{P}^{1}$.

Our main result allows us to deduce the following

\begin{corollary} \label{cor.main}
Let $V$ be an $n$-dimensional vector space over $k$.  Then there is an equivalence
$$
\mathbb{P}^{1}_{n} \equiv \mathbb{P}^{nc}(V).
$$
\end{corollary}
This corollary explains, in some sense, $\mathbb{P}^{1}_{n}$'s reliance only on $n$ and $k$ but not on the exact form of the relation $b$ in (\ref{eqn.algebraa}).

As an immediate consequence of a corollary to Theorem \ref{theorem.ft} (Corollary \ref{cor.central}) we also deduce the following related result:

\begin{corollary}
The category ${\sf C}$ has a sequence $\underline{\mathcal{L}}$ such that the six conditions above are satisfied and $D_{i}=k$ for all $i$ if and only if there is a $k$-linear equivalence ${\sf C} \equiv \mathbb{P}^{1}_{n}$ for some $n$.
\end{corollary}

\subsection{Organization of the paper}
We now briefly describe the contents of this paper.  In Section \ref{section.prelim}, after recalling the definition of $\mathbb{Z}$-algebras and noncommutative spaces of the form ${\sf Proj }A$ for $A$ a $\mathbb{Z}$-algebra, we review basic results about graded coherence and ampleness from \cite{polish}.  In Section \ref{section.ncsym}, we recall the definition of noncommutative symmetric algebras and show, in Section \ref{section.euler}, that they satisfy an Euler exact sequence, generalizing \cite[Section 3.3]{newchan} and \cite[Lemma 3.7 and Proposition 3.8]{tsen}.  We then state our main result, Theorem \ref{theorem.supermain}, in Section \ref{section.props}.  In Section \ref{section.abelian}, we study so-called linear sequences, in order to prove Theorem \ref{theorem.supermain}(1)-(4).  The argument is a refinement of the argument used to prove Theorem \ref{theorem.tsen}.

Most of the rest of the paper consists of the proof of Theorem \ref{theorem.supermain}(5), which entails proving that noncommutative projective lines are homologically well-behaved.  In order to do this, we adapt the definition and study of relative local cohomology from \cite{duality} to our context (in Section \ref{section.shom} and Section \ref{section.gor}), in order to show that noncommutative symmetric algebras are Gorenstein.  There are two fundamental differences between the analysis in \cite{duality} and this paper: since our base is affine, the relevant functors we study are much easier to work with.  On the other hand, since the noncommutative projective lines we study here are not noetherian, some of our proofs involve subtleties not encountered in the noetherian setting of \cite{duality}.

Finally, in Section \ref{section.piont}, we check that the spaces $\mathbb{P}_{n}^{1}$ satisfy the hypotheses for our main result, allowing us to deduce Corollary \ref{cor.main}.
\newline
\newline
{\it Acknowledgement:} I am grateful to the anonymous referee for their many helpful and detailed suggestions for the improvement of this paper.  Among other things, I have included the referee's proof of their significantly strengthened version of Lemma \ref{lemma.faithful}(2) and (3).

\section{Preliminaries} \label{section.prelim}

\subsection{$\mathbb{Z}$-algebras and ${\sf Proj }A$} \label{section.proja}
We recall the notion of a positively graded $\mathbb{Z}$-algebra and its graded modules, following \cite[Section 2]{quadrics}.

A $\mathbb{Z}$-{\it algebra} is a ring $A$ with decomposition $A=\oplus_{i, j \in \mathbb{Z}}A_{ij}$ into $k$-vector spaces, such that multiplication has the property $A_{ij}A_{jk} \subset A_{ik}$ while $A_{ij}A_{kl}=0$ if $j \neq k$.  Furthermore, for $i \in \mathbb{Z}$, there is a local unit $e_{i} \in A_{ii}$, such that if $a \in A_{ij}$, then $e_{i}a=a=ae_{j}$.  $A$ is {\it positively graded} if $A_{ij}=0$ for all $i>j$.  {\it In what follows, we will abuse terminology by saying '$A$ is a $\mathbb{Z}$-algebra' if $A$ is a positively graded $\mathbb{Z}$-algebra.}

If $A$ is a $\mathbb{Z}$-algebra then a graded right $A$-module $M$ is a right $A$-module together with a decomposition $M = \oplus_{i \in \mathbb{Z}}M_{i}$ such that $M_{i} A_{ij} \subset M_{j}$, $e_{i}$ acts as a unit on $M_{i}$ and $M_{i}A_{jk}=0$ if $i \neq j$.

We let ${\sf Gr }A$ denote the category of graded right $A$-modules (with obvious notion of homomorphism), and we note that it is a Grothendieck category (see \cite[Section 2]{quadrics}).

\begin{definition}
We say a positively graded $\mathbb{Z}$-algebra $A$ is a {\it connected $\mathbb{Z}$-algebra, finitely generated in degree one} if
\begin{itemize}
\item{} for all $i$, $A_{ii}$ is a division ring over $k$ for each $i$,

\item{} for all $i$, $A_{i,i+1}$ is finite-dimensional over both $A_{ii}$ and $A_{i+1,i+1}$, and

\item{} $A$ is generated in degree one, i.e. for all $i$ and for $j \geq i$, the multiplication maps $A_{ij} \otimes A_{j,j+1} \rightarrow A_{i,j+1}$ are surjective, and for $j \geq i+1$, the multiplication maps $A_{i,i+1} \otimes A_{i+1,j} \rightarrow A_{ij}$ are surjective.
\end{itemize}
\end{definition}
We remark that either part of the third item in the definition implies the other by associativity of multiplication in $A$.

Suppose $A$ is a connected $\mathbb{Z}$-algebra, finitely generated in degree one.  A graded right $A$-module $M$ is {\it right bounded} if $M_{n}=0$ for all $n>>0$.  We let ${\sf Tors }A$ denote the full subcategory of ${\sf Gr }A$ consisting of modules whose elements $m$ have the property the right $A$-module generated by $m$ is right bounded.  Then the assumption on $A$ implies that  ${\sf Tors }A$ is a Serre subcategory of ${\sf Gr A}$, and there is a torsion functor $\tau: {\sf Gr }A \rightarrow {\sf Tors }A$ which sends a module to its largest torsion submodule.  Furthermore, since ${\sf Gr }A$ has enough injectives, it follows that if $\pi:{\sf Gr }A \rightarrow {\sf Gr }A/{\sf Tors }A = : {\sf Proj }A$ is the quotient functor, then there exists a section functor $\omega:{\sf Proj }A \rightarrow {\sf Gr }A$ which is right adjoint to $\pi$.

\subsection{Graded coherence}
We now review the basic facts about coherence from \cite{polish} which we will need in the sequel.  For the rest of Section \ref{section.prelim}, we let $A$ denote a $\mathbb{Z}$-algebra such that
\begin{itemize}
\item{} $A_{ii}$ is a division ring over $k$ for each $i$, and

\item{} $A_{ij}$ is finite-dimensional as a left $A_{ii}$-module and as a right $A_{jj}$-module.
\end{itemize}
We warn the reader that our grading convention and assumptions on $A$ differ slightly from those appearing in \cite{polish}.

We let $P_{i} := \oplus_{j}A_{ij} = e_{i}A$.  We say that $M \in {\sf Gr }A$ is {\it finitely generated} if there is a surjection $P \rightarrow M$ where $P$ is a finite direct sum of modules of the form $P_{i}$.  We say $M$ is {\it coherent} if it is finitely generated and if for every homomorphism $f: P \rightarrow M$ with $P$ a finite direct sum of $P_{i}$'s, $\operatorname{ker }f$ is finitely generated.  We denote the full subcategory of ${\sf Gr }A$ consisting of coherent modules by ${\sf coh }A$.  By \cite[Proposition 1.1]{polish}, ${\sf coh }A$ is an abelian subcategory of ${\sf Gr }A$ closed under extensions.

We call $A$ {\it graded right-coherent} (or just {\it coherent}) if the right modules $P_{j}$ and $S_{j} := P_{j}/{P_{j}}_{>j}$ are coherent.  In the case that $A$ is a connected $\mathbb{Z}$-algebra, finitely generated in degree one, then the condition that the modules $P_{j}$ are coherent implies that $A$ is coherent.

We let ${\sf tors }A$ denote the full subcategory of ${\sf coh }A$ consisting of right-bounded modules.  One can check that this is a Serre subcategory of ${\sf coh }A$.  If $A$ is graded coherent, we define
$$
{\sf cohproj }A := {\sf coh }A/{\sf tors }A,
$$
which is abelian.

In the following result, which will be used in Section \ref{section.gor}, we abuse notation by letting $\pi:{\sf coh }A \rightarrow {\sf cohproj }A$ denote the quotient functor.
\begin{lemma} \label{lemma.faithful}
Suppose $A$ is a coherent connected $\mathbb{Z}$-algebra, finitely generated in degree one. Then
\begin{enumerate}
\item{} the inclusion functor $\iota:{\sf coh }A \rightarrow {\sf Gr }A$ descends to an exact functor
$$
\underline{\iota}:{\sf cohproj }A \rightarrow {\sf Proj }A,
$$

\item{} the functor $\underline{\iota}$ is fully faithful,

\item{} if $M$ and $N$ are objects in ${\sf coh }A$ such that the torsion submodules of $\iota(M)$ and $\iota(N)$ are coherent, then the map
\begin{equation} \label{eqn.extinducee}
\operatorname{Ext}^{1}_{{\sf cohproj }A}(\pi(M), \pi(N)) \rightarrow \operatorname{Ext}^{1}_{{\sf Proj }A}(\underline{\iota}(\pi(M)), \underline{\iota}(\pi(N)))
\end{equation}
induced by $\underline{\iota}$ is an isomorphism of $\operatorname{End }(\pi(N))-\operatorname{End }(\pi(M))$-bimodules.
\end{enumerate}
\end{lemma}

\begin{proof}
We first note that, by hypothesis, ${\sf Proj }A$ is well-defined.

By \cite[Proposition 1.1]{polish}, the inclusion $\iota: {\sf coh }A \rightarrow {\sf Gr }A$ makes ${\sf coh }A$ an abelian subcategory of ${\sf Gr }A$.  Furthermore, if $M \in {\sf tors }A$, then $\iota(M) \in {\sf Tors }A$.  Therefore, by \cite[Corollaire 2, p. 368]{gab}, $\iota$ induces a functor $\underline{\iota}:{\sf cohproj }A \rightarrow {\sf Proj }A$ which is exact by \cite[Corollaire 3, p. 369]{gab}.  Part (1) of the lemma follows.

We next show that $\underline{\iota}$ is faithful.  To this end, we recall that, by definition of $\underline{\iota}$, $\underline{\iota}(\pi(M)) = \pi(\iota(M))$, and if $f \in \operatorname{Hom }_{{\sf cohproj }A}(\pi(M), \pi(N))$, then we may identify $f$ with a map $f' \in \operatorname{Hom }_{{\sf coh }A}(M', N/N')$ where $M' \subset M$ is coherent, $M/M' \in {\sf tors }A$, and $N' \subset N$ is in ${\sf tors }A$.  Furthermore, by the proof of \cite[Corollaire 2, p. 368]{gab}, $\underline{\iota}(f)\in \operatorname{Hom }_{{\sf Proj }A}(\pi(\iota(M)), \pi(\iota(N)))$ corresponds to $\iota(f') \in \operatorname{Hom }_{{\sf Gr }A}(\iota(M'), \iota(N)/\iota(N'))$.  Therefore, if $\underline{\iota}(f)=0$, it follows that $\pi(\iota(f'))=0$ which implies that the image of $\iota(f')$ is torsion.  Since $M'$ is coherent, the image of $\iota(f')$ is right-bounded, so that the image of $f'$ is in ${\sf tors }A$.  It follows that $\pi(f')=0$, so that $f=0$.

Now we prove $\underline{\iota}$ is full.  Suppose an element of $\operatorname{Hom }_{{\sf Proj }A}(\underline{\iota} \pi(M), \underline{\iota} \pi (N))$ corresponds to a map $g:M' \rightarrow \iota(N)/N'$ in ${\sf Gr }A$, where
$N'$ and $\iota(M)/M'$ are in ${\sf Tors} A$.  We first claim that $\iota(M)$ and $M'$  are equal in high degree and that $M'$ is coherent.  For, since $M$ is coherent, it is finitely generated, so that $\iota(M)$ and $M'$ are equal in high degree.  Thus, by our assumptions on $A$, $\iota(M)/M'$ is coherent.  It follows that $M'$ is coherent, whence the claim.  Thus, without loss of generality, we may assume $\iota(M)=M'$.  Now suppose $P$ is a finite direct sum of modules of the form $P_{i}$ and $K$ is a coherent submodule of $P$ such that $\iota(M) \cong P/K$.  Since $P$ is projective, $g$ lifts to a map $g':P \rightarrow \iota(N)$.  Furthermore, since $g'(K)$ is the image of a map $K \rightarrow N$ of coherent modules, $g'(K)$ is coherent. It follows that $g$ factors through a map $\iota(M) \rightarrow \iota(N)/g'(K)$, so comes from a morphism in ${\sf cohproj }A$.

Finally, we prove (3).  We first recall that the group $\operatorname{Ext}^{1}_{{\sf cohproj }A}(\pi(M), \pi(N))$ is defined as Yoneda's extension group, while the group $\operatorname{Ext}^{1}_{{\sf Proj }A}(\underline{\iota}(\pi(M)), \underline{\iota}(\pi(N)))$ defined in terms of injective resolutions is isomorphic to Yoneda's extension group as a bimodule \cite[Section VII.7]{mitchell}.

By (1), the functor $\underline{\iota}$ applied to an exact sequence in ${\sf cohproj }A$ is exact in ${\sf Proj }A$, and hence determines an element of $\operatorname{Ext}^{1}_{{\sf Proj }A}(\underline{\iota}(\pi(M)), \underline{\iota}(\pi(N)))$.  The fact that this induces a morphism between extension groups follows immediately from functoriality of $\underline{\iota}$.

Since $\underline{\iota}$ is exact and additive it preserves pullbacks and pushouts.  It follows that $\underline{\iota}$ induces an additive function between extension groups, and that the assignment is compatible with bimodule structures, as one can check.

Next we note that by hypothesis we may assume, without loss of generality, that $\iota(M)$ and $\iota(N)$ are torsion-free.  To prove injectivity of (\ref{eqn.extinducee}), suppose $\underline{\iota}$ applied to an extension
$$
0 \rightarrow \pi(N) \rightarrow \pi(P) \rightarrow \pi(M) \rightarrow 0
$$
in ${\sf cohproj }A$ maps to a trivial extension in ${\sf Proj }A$.  We then have an isomorphism $\underline{\iota}\pi (P) \rightarrow \underline{\iota} \pi (M \oplus N)$ in ${\sf Proj }A$.  By part (2), the isomorphism must be induced by an isomorphism $\pi(P) \rightarrow \pi (M \oplus N)$, and the result follows.

Next, we prove that (\ref{eqn.extinducee}) is surjective.  In order to simplify the exposition, we omit the functors $\iota$ and $\underline{\iota}$.  Suppose
\begin{equation} \label{eqn.origext}
0 \rightarrow \pi(N) \rightarrow \pi(E) \overset{f}{\rightarrow} \pi(M) \rightarrow 0
\end{equation}
represents an element of $\operatorname{Ext}^{1}_{{\sf Proj }A}(\pi(M), \pi(N))$, and suppose
\begin{equation} \label{eqn.mpres}
0 \rightarrow K_{M} \rightarrow P \rightarrow M \rightarrow 0
\end{equation}
is a finite presentation of $M$.  Thus, $P$ is a finite direct sum of modules of the form $P_{i}$ so that $P$ and $K_{M}$ are coherent.  We consider the commutative diagram in ${\sf Gr }A$
\begin{equation} \label{eqn.doubled}
\begin{CD}
0 & \longrightarrow & N & \longrightarrow & P \oplus N & \longrightarrow & P & \longrightarrow & 0 \\
& & @VVV @V{\phi}VV @VV{\phi_{M}}V \\
0 & \longrightarrow & \omega \pi(N) & \longrightarrow & \omega \pi(E) & \longrightarrow & \mbox{im }\omega f & \longrightarrow & 0
\end{CD}
\end{equation}
where the top row is the canonical split exact sequence, the bottom row is induced by $\omega$ applied to (\ref{eqn.origext}), the left vertical is the adjointness map, the right vertical is induced by the composition of the right arrow in (\ref{eqn.mpres}) with the adjointness map $M \rightarrow \omega \pi M$, and the center vertical is due to the projectivity of $P$.  We begin by showing that $\mbox{im }\phi$ is coherent.  Let $K =\operatorname{ker }\phi$.  From the kernel-cokernel exact sequence of (\ref{eqn.doubled}) we deduce that $K_{M}/K$ is torsion.  As in the proof of (2), the fact that $K_{M}$ is coherent implies that $K$ is coherent.  Thus, $\mbox{im }\phi \cong P \oplus N/K$ is coherent as well.

To complete the proof of (3), we claim that the original extension (\ref{eqn.origext}) is equivalent to an extension of the form
$$
0 \rightarrow \pi(N) \rightarrow \pi (\operatorname{im }\phi) \rightarrow \pi(M) \rightarrow 0
$$
in ${\sf Proj }A$, so that surjectivity of (\ref{eqn.extinducee}) will follow from (2).  To prove the claim, let $\psi:\operatorname{im }\phi \rightarrow M$ denote the restriction of the bottom right horizontal in (\ref{eqn.doubled}) and consider the diagram
$$
\begin{CD}
0 & \longrightarrow & \pi(\operatorname{ker }\psi) & \longrightarrow & \pi(\operatorname{im }\phi) & \longrightarrow & \pi(M) & \longrightarrow 0 \\
& & @VVV @VVV @VVV \\
0 & \longrightarrow & \pi \omega \pi(N) & \longrightarrow & \pi \omega \pi(E) & \longrightarrow & \pi \omega \pi(M) & \longrightarrow & 0
\end{CD}
$$
in which the bottom row is induced by (\ref{eqn.origext}) and the verticals are induced by inclusion.  Then the left and right verticals are isomorphisms since $N \subset \operatorname{ker }\psi \subset \omega \pi(N)$ and the kernel and cokernel of the adjointness map is torsion.  It follows that the center vertical is an isomorphism and the result follows in light of the fact that $\pi \omega$ is isomorphic to the identity functor.
\end{proof}

\subsection{Ampleness} \label{section.ampleness}
We let ${\sf C}$ denote a $k$-linear category, we let $\underline{\mathcal{E}} = (\mathcal{E}_{i})_{i \in \mathbb{Z}}$ denote a sequence of objects in ${\sf C}$ such that $\operatorname{Hom}(\mathcal{E}_{i},\mathcal{E}_{i})$ is a division ring, and (in this section) we assume that for all $\mathcal{M} \in {\sf C}$, the dimension of $\operatorname{Hom}_{\sf C}(\mathcal{E}_{i},\mathcal{M})$ is finite as a right $\operatorname{End}(\mathcal{E}_{i})$-module.

We call $\underline{\mathcal{E}}$
\begin{itemize}
\item{} {\it projective} if for every surjection $f: \mathcal{M} \rightarrow \mathcal{N}$ in ${\sf C}$ there exists an integer $n$ such that $\operatorname{Hom}_{\sf C}(\mathcal{E}_{-i},f)$ is surjective for all $i>n$, and

\item{} {\it ample} if it is projective, and if for every $\mathcal{M} \in {\sf C}$ and every $m \in \mathbb{Z}$ there exists a surjection
    $$
    \oplus_{j=1}^{s}\mathcal{E}_{-i_{j}} \rightarrow \mathcal{M}
    $$
for some $i_{1},\ldots, i_{s}$ with $i_{j} \geq m$ for all $j$.
\end{itemize}

We will need the following mild generalization of (part of) \cite[Proposition 2.3(ii) and Theorem 2.4]{polish}, which holds by the same proof:

\begin{theorem} \label{theorem.polish}
If $\underline{\mathcal{E}}$ is an ample sequence and $A(\underline{\mathcal{E}})$ denotes the $\mathbb{Z}$-algebra with $A_{ij}:=\operatorname{Hom }_{\sf C}(\mathcal{E}_{-j}, \mathcal{E}_{-i})$ for $i \leq j$ and $A_{ij}=0$ otherwise, then $A(\underline{\mathcal{E}})$ is coherent, and there is an equivalence
$$
{\sf C} \equiv {\sf cohproj }A(\underline{\mathcal{E}}).
$$
\end{theorem}

\section{Noncommutative symmetric algebras} \label{section.ncsym}
Let $B$ and $C$ be noetherian $k$-algebras.  In this section, following \cite{vandenbergh}, we define the noncommutative symmetric algebra of certain $B-C$-bimodules.  Some of the exposition is adapted from \cite{witt}.

\subsection{Bimodules}
We assume throughout this section that $N$ is an $B-C$-bimodule which is finitely generated projective as both a left $B$-module and as a right $C$-module.  We recall that the {\it right dual of $N$}, denoted $N^{*}$, is the $C-B$-bimodule whose underlying set is
$\operatorname{Hom}_{C}(N_{C},C)$, with action
$$
(c \cdot \psi \cdot b)(n)=c\psi(bn)
$$
for all $\psi \in \operatorname{Hom}_{C}(N_{C},C)$, $c\in C$ and $b \in B$.

The {\it left dual of $N$}, denoted ${}^{*}N$, is the $C-B$-bimodule whose underlying set is
$\operatorname{Hom}_{B}({}_{B}N,B)$, with action
$$
(c \cdot \phi
\cdot b)(n)=\phi(nc)b
$$
for all $\phi \in \operatorname{Hom}_{B}({}_{B}N,B)$, $c \in C$ and $b \in B$.  This assignment extends to morphisms between $B-C$-bimodules in the obvious way.

We set
$$
N^{i*}:=
\begin{cases}
N & \text{if $i=0$}, \\
(N^{i-1*})^{*} & \text{ if $i>0$}, \\
{}^{*}(N^{i+1*}) & \text{ if $i<0$}.
\end{cases}
$$
In general, $N$ may not be isomorphic to $N^{**}$ or ${}^{**}N$ \cite[Section 6.4]{ringel}.  Furthermore, although $N^{*}$ (resp. ${}^{*}N$) is finitely generated projective on the left (resp. finitely generated projective on the right), it is not clear that $N^{*}$ is finitely generated projective on the right (resp. finitely generated projective on the left).  Therefore, we make the following

\begin{definition}
We say $N$ is {\it admissible} if $N^{i*}$ is finitely generated projective on each side for all $i \in \mathbb{Z}$.  We say $N$ is {\it 2-periodic} if $N$ is admissible, $N^{i*}$ is free on each side, and the left rank of $N^{i*}$ equals the right rank of $N^{i+1*}$.
\end{definition}
We remark that if $B$ and $C$ are finite dimensional simple rings over $k$, then $N$ is automatically 2-periodic.  If $B$ and $C$ are fields and $N$ is of type $(2,2)$, then $N$ is 2-periodic \cite[Lemma 3.4]{newchan}.  Finally, if $B=C$ is a perfect field and $N$ has finite left and right dimension, then $N$ is 2-periodic \cite[Proposition 4.3]{hart}.

If $B$ and $C$ are division rings and $N$ has finite left-dimension $m$ and finite right-dimension $n$, we say $N$ has {\it left-right dimension }$(m,n)$.  In this situation, we let $\operatorname{ldim}N$ denote the dimension of $N$ over $B$ and we let $\operatorname{rdim}N$ denote the dimension of $N$ over $C$.

Let
$$
S_{i} = \begin{cases}
B & \mbox{if $i$ is even, and} \\
C & \mbox{if $i$ is odd.}
\end{cases}
$$
In what follows, all unadorned tensor products will be over $S_{i}$.

We recall that, if $N$ is admissible, then for each $i$, both pairs of functors
\begin{equation} \label{eqn.adjointone}
(-\otimes_{S_{i}} N^{i*},-\otimes_{S_{i+1}} N^{i+1*})
\end{equation}
and
\begin{equation} \label{eqn.adjointtwo}
(-\otimes_{S_{i}} {}^{*}(N^{i+1*}), -\otimes_{S_{i+1}} N^{i+1*})
\end{equation}
between the category of right $S_{i}$-modules and the category of right $S_{i+1}$-modules have canonical adjoint structures.

By adjointness, $S_{i}$ maps to $N^{i*} \otimes_{S_{i+1}} N^{i+1*}$ and to ${}^{*}(N^{i+1*}) \otimes_{S_{i+1}} N^{i+1*}$ and we denote its images by $Q_{i}$ and $Q_{i}'$, respectively.  If $B$ and $C$ are division rings, $N$ is 2-periodic, $\{\phi_{1}, \ldots, \phi_{n}\}$ is a right-basis for $N^{i*}$ and $\{\phi_{1}^{*}, \ldots, \phi_{n}^{*}\}$ is a corresponding left dual basis for $N^{i+1*}$, then the canonical adjointness map from $S_{i}$ to $Q_{i}$ maps $1$ to $\sum_{i}\phi_{i} \otimes \phi_{i}^{*}$ (see \cite[Section 2]{tsen}).  In particular, the latter element is $S_{i}$-central.  We will employ this fact without comment in the sequel.

\subsection{The definition of $\mathbb{S}^{nc}(M)$}
 We now recall (from \cite{vandenbergh}) the definition of the noncommutative symmetric algebra of an admissible $B$-$C$-bimodule $N$.  The {\it noncommutative symmetric algebra of $N$}, denoted $\mathbb{S}^{nc}(N)$, is the positively graded $\mathbb{Z}$-algebra $\underset{i, j \in \mathbb{Z}}{\oplus}A_{ij}$ with components defined as follows:
\begin{itemize}
\item{} $A_{ii}=B$ for $i$ even,

\item{} $A_{ii}=C$ for $i$ odd, and

\item{} $A_{i,i+1}=N^{i*}$.
\end{itemize}
In order to define $A_{ij}$ for $j>i+1$, we introduce some notation: we define $T_{i,i+1} := A_{i,i+1}$, and, for $j>i+1$, we define
$$
T_{ij} := A_{i,i+1} \otimes A_{i+1,i+2} \otimes \cdots \otimes A_{j-1,j}.
$$
We let $R_{i,i+1}:= 0$, $R_{i,i+2}:=Q_{i}$,
$$
R_{i,i+3}:=Q_{i} \otimes N^{i+2*}+N^{i*} \otimes Q_{i+1},
$$
and, for $j>i+3$, we let
$$
R_{ij} := Q_{i} \otimes T_{i+2,j}+T_{i,i+1}\otimes Q_{i+1} \otimes T_{i+3,j}+\cdots + T_{i,j-2} \otimes Q_{j-2}.
$$

\begin{itemize}
\item{} For $j>i+1$, we define $A_{ij}$ as the quotient $T_{ij}/R_{ij}$.
\end{itemize}
Multiplication in $\mathbb{S}^{nc}(N)$ is defined as follows:
\begin{itemize}
\item{} if $x \in A_{ij}$ and $y \in A_{jk}$, with either $i=j$ or $j=k$, then $xy$ is induced by the usual scalar action,

\item{}  otherwise, if $i<j<k$, we have
\begin{eqnarray*}
A_{ij} \otimes A_{jk} & = & \frac{T_{ij}}{R_{ij}} \otimes \frac{T_{jk}}{R_{jk}}\\
& \cong & \frac{T_{ik}}{R_{ij}\otimes T_{jk}+T_{ij} \otimes
R_{jk}}.
\end{eqnarray*}
Since $R_{ij} \otimes T_{jk}+T_{ij} \otimes R_{jk}$ is a submodule of $R_{ik}$,
we may define multiplication $A_{ij} \otimes A_{jk} \longrightarrow
A_{ik}$ as the canonical epimorphism.
\end{itemize}

\subsection{Euler sequences} \label{section.euler}
In this section, we assume $D_{0}$ and $D_{1}$ are division rings over $k$, we let $M$ be a $k$-central $D_{0}-D_{1}$-bimodule such that $M$ is 2-periodic and not of type $(1,1)$, $(1,2)$ or $(1,3)$, and we let $A=\mathbb{S}^{nc}(M)$.  Our main goal in this section, Corollary \ref{cor.euler}, is to prove that the trivial module in ${\sf Gr }A$ has the expected resolution.  We then derive some consequences, which we will need in the sequel.  We remark that if $M$ is 2-periodic, of type $(1,1)$, $(1,2)$ or $(1,3)$, and Corollary \ref{cor.euler} holds for $A$, then it is not hard to show that $A$ degenerates so that ${\sf cohproj }A$ is trivial.

\begin{lemma} \label{lemma.zerodivisor}
Let $j \in \mathbb{Z}$ and suppose $M^{j*}$ has left-right dimension $(n,m)$ with $m \geq 2$.  If $v \in M^{j*}$ has the property that $v \otimes g \in Q_{j}$ for some nonzero $g \in M^{j+1*}$, then $v =0$.
\end{lemma}

\begin{proof}
Let $\{\phi_{1}, \ldots, \phi_{m}\}$ denote a right basis for $M^{j*}$ so that $\{\phi_{1}^{*}, \ldots, \phi_{m}^{*}\}$ (the right-dual basis) is a left basis for $M^{j+1*}$.  Let $v=\sum_{l}\phi_{l}a_{l}$ and suppose $v \neq 0$.  We have
\begin{eqnarray*}
(\sum_{l} \phi_{l} \otimes \phi_{l}^{*})c & = & \sum_{l} \phi_{l} a_{l} \otimes g \\
& = & \sum_{l} \phi_{l} \otimes a_{l}g.
\end{eqnarray*}
Therefore, for all $l$, $a_{l}g = \phi_{l}^{*}c$.  Since there exists an $l$ such that $a_{l} \neq 0$, we have $c \neq 0$.  It follows that all $a_{l}$ are nonzero.  Since $m \geq 2$, we thus have $a_{1}^{-1}\phi_{1}^{*}=a_{2}^{-1}\phi_{2}^{*}$, which is a contradiction.
\end{proof}
The following generalizes \cite[Lemma 3.5]{newchan}:

\begin{prop} \label{prop.zerodivisor}
Let $j \in \mathbb{Z}$ and suppose $M^{j*}$ has left-right dimension $(n,m)$ with $m,n \geq 2$.  If $v \in T_{i,j+1}$ has the property that $v \otimes g \in R_{i,j+2}$ for some nonzero $g \in M^{j+1*}$, then $v \in R_{i,j+1}$.
\end{prop}

\begin{proof}
Throughout the proof, we let $\{f_{1},\ldots,f_{n}\}$ be a left basis of $M^{j*}$ so that $\{{}^{*}f_{1}, \ldots, {}^{*}f_{n}\}$ (the left-dual basis) is a right basis of $M^{j-1*} \cong {}^{*}(M^{j*})$.  Similarly, we let $\{\phi_{1}, \ldots, \phi_{m}\}$ denote a right basis for $M^{j*}$ so that $\{\phi_{1}^{*}, \ldots, \phi_{m}^{*}\}$ (the right-dual basis) is a left basis for $M^{j+1*}$.

We proceed by induction on $j-i$.  First, if $i=j$, then the result follows from Lemma \ref{lemma.zerodivisor}.

Next, we assume $j=i+1$.  If $g=\sum_{r}b_{r}\phi_{r}^{*}$, then there exists an $h \in M^{j-1*}$ such that
$$
v \otimes \sum_{r}b_{r}\phi_{r}^{*} - h \otimes \sum_{s} \phi_{s} \otimes \phi_{s}^{*} \in R_{i,j+1} \otimes M^{j+1*}.
$$
It follows that
$$
vb_{r}-h \otimes \phi_{r} \in R_{i,j+1}
$$
for all $r$.  Thus
$$
vb_{r}-h \otimes \phi_{r} = c_{r}(\sum_{l}{}^{*}f_{l} \otimes f_{l})
$$
for all $r$.  Without loss of generality, assume $b_{1} \neq 0$.  Then we deduce
$$
h \otimes \phi_{1}b_{1}^{-1}b_{2}+c_{1}b_{1}^{-1}b_{2}(\sum_{l}{}^{*}f_{l} \otimes f_{l})=h \otimes \phi_{2}+c_{2}(\sum_{l}{}^{*}f_{l} \otimes f_{l}).
$$
Thus,
$$
h \otimes (\phi_{1}b_{1}^{-1}b_{2}-\phi_{2})+(c_{1}b_{1}^{-1}b_{2}-c_{2})(\sum_{l}{}^{*}f_{l}\otimes f_{l})=0.
$$
Since the $\phi_{i}$ are right-independent, $\psi:= \phi_{1}b_{1}^{-1}b_{2}-\phi_{2} \neq 0$.  Therefore, $h \otimes \psi \in R_{i,i+2}$, so that, by Lemma \ref{lemma.zerodivisor}, $h=0$.  This implies that $v \in R_{i,j+1}$ as desired.
%If we write $\psi=\sum d_{q}f_{q}$, then we have
%$$
%h \otimes \psi = h \otimes \sum_{q} d_{q}f_{q} = e \sum {}^{*}f_{l} \otimes f_{l}
%$$
%and, therefore, $hd_{q}=e{}^{*}f_{q}$ for all $q$.  If $d_{q}=0$ for some $q$, then $e=0$ so that $hd_{q}=0$ for all $q$, which is a contradiction unless $h=0$.  Otherwise, $d_{1} \neq 0$ and $d_{2} \neq 0$, which implies that since $hd_{1}=e{}^{*}f_{1}$ and $hd_{2}=e{}^{*}f_{2}$, we have $e{}^{*}f_{1}d_{1}^{-1}=e{}^{*}f_{2}d_{2}^{-1}$.  If $e \neq 0$, this contradicts the right-independence of the ${}^{*}f_{i}$'s.  Therefore, $e=0$ so that $h=0$ which implies that $v \in R_{i j+1}$ as desired.

Now we assume $j>i+1$ and let $g=\sum_{r}b_{r}\phi_{r}^{*}$ as above.  Then
$$
v \otimes g - h \otimes \sum \phi_{s} \otimes \phi_{s}^{*} \in R_{i,j+1} \otimes M^{j+1*}
$$
for some $h \in T_{i j}$ so that $vb_{r}-h \otimes \phi_{r} \in R_{i,j+1}$ for all $r$.  Thus, for all $r$, we have
$$
vb_{r}-h \otimes \phi_{r} - h_{r} \otimes (\sum {}^{*}f_{l} \otimes f_{l}) \in R_{ij} \otimes M^{j*}
$$
for some $h_{r} \in T_{i,j-1}$.  Suppose, without loss of generality, that $b_{1} \neq 0$.  Then
\begin{equation} \label{eqn.long1}
vb_{1}-h \otimes \phi_{1}-h_{1} \otimes (\sum {}^{*}f_{l} \otimes f_{l}) \in R_{ij} \otimes M^{j*}
\end{equation}
while
\begin{equation} \label{eqn.long2}
vb_{2}-h \otimes \phi_{2} - h_{2} \otimes ({}^{*}f_{l} \otimes f_{l}) \in R_{ij} \otimes M^{j*}.
\end{equation}
Therefore, if we multiply (\ref{eqn.long1}) by $b_{1}^{-1}b_{2}$ on the right and subtract (\ref{eqn.long2}), we deduce that
$$
h \otimes (\phi_{2}-\phi_{1}b_{1}^{-1}b_{2})+(h_{2}-h_{1}b_{1}^{-1}b_{2})\otimes (\sum {}^{*}f_{l} \otimes f_{l}) \in R_{ij} \otimes M^{j*}.
$$
Since $\phi_{2}-\phi_{1}b_{1}^{-1}b_{2} \neq 0$, induction implies that $h \in R_{ij}$, and it follows that $v \in R_{i,j+1}$.
\end{proof}

\begin{theorem} \label{theorem.euler}
For all $i \leq j$, the canonical complex
\begin{equation} \label{eqn.seqq}
0 \rightarrow A_{ij} \otimes Q_{j} \rightarrow A_{i,j+1} \otimes M^{j+1*} \rightarrow A_{i,j+2} \rightarrow 0
\end{equation}
is exact.
\end{theorem}

\begin{proof}
There are two fundamentally different cases to consider.  First, if $M$ has left-right dimension $(1,n)$ for $n \geq 4$, then the proofs of \cite[Lemma 3.7]{tsen} and \cite[Proposition 3.8]{tsen} still work in our context.  Thus, we may suppose that $M$ is of type $(n,m)$ with both $n$ and $m$ $\geq 2$.  In this case, suppose without loss of generality, that $M^{j*}$ has left-right dimension $(n,m)$.  We show that (\ref{eqn.seqq}) is exact on the left.  In order to prove this, it suffices to check that
\begin{equation} \label{eqn.int}
R_{i,j+1} \otimes M^{j+1*} \cap T_{ij} \otimes Q_{j} = R_{ij} \otimes Q_{j}.
\end{equation}
As in the proof of Proposition \ref{prop.zerodivisor}, we let $\{f_{1},\ldots,f_{n}\}$ be a left basis of $M^{j*}$ so that $\{{}^{*}f_{1}, \ldots, {}^{*}f_{n}\}$ (the left-dual basis) is a right basis of $M^{j-1*} \cong {}^{*}(M^{j*})$.  Similarly, we let $\{\phi_{1}, \ldots, \phi_{m}\}$ denote a right basis for $M^{j*}$ so that $\{\phi_{1}^{*}, \ldots, \phi_{m}^{*}\}$ (the right-dual basis) is a left basis for $M^{j+1*}$.

First, we assume $j=i$.  Then both sides of (\ref{eqn.int}) are zero, as desired.

Next, we assume $j=i+1$.  Suppose $v \otimes \sum_{p} \phi_{p} \otimes \phi_{p}^{*}$ is an element of the left-hand side of (\ref{eqn.int}).  Then we have an equality
$$
v \otimes \sum_{p} \phi_{p} \otimes \phi_{p}^{*} = a (\sum_{q} {}^{*}f_{q} \otimes f_{q}) \otimes \sum_{r} a_{r} \phi_{r}^{*}
$$
for some scalars $a, a_{1}, \ldots, a_{m}$, so that $v \otimes \phi_{l} = a a_{l} \sum_{q} {}^{*}f_{q} \otimes f_{q}$ for all $l$.  Without loss of generality, we may assume $\phi_{1}=f_{1}$ so that $a a_{1}=0$.  If $a=0$ then $v \otimes \phi_{l}=0$ for all $l$ which implies that $v=0$.  Otherwise, $a_{1}=0$, so that $v \otimes f_{1}=0$.  By Proposition \ref{prop.zerodivisor}, $v=0$ and the assertion follows in this case.

Now suppose $j>i+1$.  If $v \otimes \sum_{p} \phi_{p} \otimes \phi_{p}^{*}$ is in the left-hand side of (\ref{eqn.int}), then it is in $R_{i, j+1} \otimes M^{j+1*}$.  If we let $\mbox{ev}_{\phi_{q}} \in {}^{*}(M^{j+1*})$ denote the evaluation map, then $\mbox{id} \otimes \mbox{ev}_{\phi_{q}}(v \otimes \sum_{p} \phi_{p} \otimes \phi_{p}^{*})=v \otimes \phi_{q} \in R_{i, j+1}$.  Therefore, by Proposition \ref{prop.zerodivisor}, $v \in R_{ij}$.

%Then there exists a $w \in R_{ij}\otimes M^{j*}$, a $u \in T_{ij-2}$ and scalars $b_{1}, \ldots, b_{m}$ such that
%$$
%v \otimes \sum_{p} \phi_{p} \otimes \phi_{p}^{*}=\sum_{q} w \otimes a_{q}\phi_{q}^{*}+u \otimes (\sum_{r}{}^{*}f_{r} \otimes f_{r}) \otimes \sum_{s} b_{s} \phi_{s}^{*}.
%$$
%It follows that, for all $1 \leq l \leq m$, we have
%\begin{equation} \label{eqn.newt}
%v \otimes \phi_{l}-u b_{l} \otimes \sum_{r} {}^{*}f_{r} \otimes f_{r} \in R_{ij} \otimes M^{j*}.
%\end{equation}
%If $b_{l}=0$ for some $l$, Proposition \ref{prop.zerodivisor} thus implies $v \in R_{ij}$ as desired.  Thus, assume that $b_{l} \neq 0$ for all $l$.  Since, without loss of generality, we may take $\phi_{1}=f_{1}$, (\ref{eqn.newt}) implies that $v -ub_{1}\otimes {}^{*}f_{1} \in R_{ij}$ and $ub_{1}\otimes {}^{*}f_{2} \in R_{ij}$.  The second inclusion and Proposition \ref{prop.zerodivisor} implies that $u \in R_{ij-1}$ so that the first inclusion implies that $v \in R_{ij}$ as desired.
\end{proof}

\begin{corollary} \label{cor.euler}
For all $k \in \mathbb{Z}$, multiplication in $A$ induces an exact sequence of $D_{i-2}-A$-bimodules
$$
0 \rightarrow Q_{i-2} \otimes e_{i}A \rightarrow A_{i-2,i-1} \otimes e_{i-1}A \rightarrow e_{i-2}A \rightarrow e_{i-2}A/e_{i-2}A_{\geq i-1} \rightarrow 0.
$$
\end{corollary}

\begin{proof}
The only nontrivial part of the proof is to show the sequence is exact on the left.  This can be checked by counting right-dimensions, which can be deduced from the exactness of (\ref{eqn.seqq}).  More precisely, one can prove, using (\ref{eqn.seqq}) and induction on the difference of indices, that the right-dimension of $A_{i-2,j}$ equals the right-dimension of $A_{i,j+2}$, and the right-dimension of $A_{i-2,i-1} \otimes A_{i-1,j}$ equals the right-dimension of $A_{i,j+1} \otimes A_{j+1,j+2}$.  Using these facts, the result follows immediately from (\ref{eqn.seqq}).  We leave the details to the reader.
\end{proof}

We will need the following result for our proof of Theorem \ref{thm.gorenstein}.

\begin{lemma} \label{lemma.euler}
If $x \in A_{i,j+1}$ is such that $x y = 0$ for all $y \in A_{j+1,j+2}$, then $x=0$.
\end{lemma}

\begin{proof}
If $j+1 \leq i$ the result is trivial, so suppose $j+1 > i$.  If $A_{j+1,j+2}$ has left-right dimension $(n,m)$ with $m,n \geq 2$, the result follows from Proposition \ref{prop.zerodivisor}.  Similarly, if $A_{j+1,j+2}$ has left-right dimension $(1,n)$ with $n \geq 4$, the result follows from \cite[Lemma 3.7]{tsen}.

It remains to prove the result in case $A_{j+1,j+2}$ has left-right dimension $(n,1)$ with $n \geq 4$.  We let $\{\phi_{1}, \ldots, \phi_{m}\}$ denote a right basis for $A_{j,j+1}$ so that $\{\phi_{1}^{*}, \ldots, \phi_{m}^{*}\}$ (the right-dual basis) is a left basis for $M^{j+1*}$.  Since $xy=0$ for all $y \in A_{j+1, j+2}$, Theorem \ref{theorem.euler} implies that in $A_{i,j+1} \otimes M^{j+1*}$, $x \otimes \phi_{1}^{*} = \sum_{l} z \phi_{l} \otimes \phi_{l}^{*}$ for some $z \in A_{ij}$.  It follows that $x=z\phi_{1}$, $0=z \phi_{2}$ and $0=z \phi_{3}$.  By \cite[Lemma 3.7]{tsen}, $z=0$, so that $x=0$.
\end{proof}

Since $A$ is a connected $\mathbb{Z}$-algebra, finitely generated in degree one, we may form the category ${\sf Proj }A$.  As in Section \ref{section.proja}, we let $\pi: {\sf Gr }A \rightarrow {\sf Proj }A$ denote the quotient functor, and we let $\mathcal{A}_{i} := \pi (e_{i}A)$.  Applying $\pi$ to the Euler sequence from Corollary \ref{cor.euler} yields, for each $i \in \mathbb{Z}$, an exact sequence
\begin{equation} \label{eqn.ei}
0 \rightarrow \mathcal{A}_{i} \rightarrow \mathcal{A}_{i-1}^{d_{i}} \rightarrow \mathcal{A}_{i-2} \rightarrow 0
\end{equation}
in ${\sf Proj }A$, where $d_{i}$ is the right-dimension of $A_{i-2,i-1}$.

\begin{prop} \label{prop.phiiso}
Suppose $a \in \operatorname{End }(\mathcal{A}_{i})$ is induced by left-multiplication by $\alpha \in A_{ii}$.  Then there exist $f_{a} \in \operatorname{End}(\mathcal{A}_{i-1}^{d_{i}})$ and $a' \in \operatorname{End }(\mathcal{A}_{i-2})$ such that the diagram
\begin{equation} \label{eqn.doublesquare}
\begin{CD}
0 & \longrightarrow & \mathcal{A}_{i} & \longrightarrow & \mathcal{A}_{i-1}^{d_{i}} & \longrightarrow & \mathcal{A}_{i-2} & \longrightarrow 0  \\
& & @V{a}VV  @VV{f_{a}}V @VV{a'}V \\
0 & \longrightarrow & \mathcal{A}_{i} & \longrightarrow & \mathcal{A}_{i-1}^{d_{i}} & \longrightarrow & \mathcal{A}_{i-2} & \longrightarrow 0,
\end{CD}
\end{equation}
whose horizontals are (\ref{eqn.ei}), commutes.  Moreover, $a'$ is induced by left-multiplication by $\alpha \in A_{i-2, i-2}=A_{i i}$.
\end{prop}

\begin{proof}
We prove the existence of the desired morphisms in the graded module category.  Since $Q_{i-2}$ is $D_{i-2}$-central, the result follows from the fact that the morphisms in the exact sequence from Corollary \ref{cor.euler} are compatible with left-multiplication by elements of $D_{i-2}$.

\end{proof}

\section{Statement of the main theorem}  \label{section.props}
In this section, ${\sf C}$ will denote a $k$-linear abelian category.  Our main theorem will involve a collection of properties describing a sequence
$$
\underline{\mathcal{L}} := (\mathcal{L}_{i})_{i \in \mathbb{Z}}
$$
of objects in ${\sf C}$.  Before we describe these properties, we recall that if $\mathcal{L}$ and $\mathcal{M}$ are objects of ${\sf C}$ and $\operatorname{End }\mathcal{L} = :D$, then one can define an object ${}^{*}\operatorname{Hom}(\mathcal{M}, \mathcal{L}) \otimes_{D} \mathcal{L}$ of ${\sf C}$ as in \cite[Section B3]{az2}.  In particular, if $D$ is a noetherian $k$-algebra and $\operatorname{Hom}(\mathcal{M}, \mathcal{L})$ is a finitely generated projective left $D$-module, then there is a canonical map
$$
\eta_{\mathcal{M}}: \mathcal{M} \rightarrow {}^{*}\operatorname{Hom}(\mathcal{M}, \mathcal{L}) \otimes_{D} \mathcal{L}.
$$
%Furthermore, if $\mathcal{N}$ is an object of ${\sf C}$ and $\operatorname{Hom}(\mathcal{M}, \mathcal{L})$ and $\operatorname{Hom}(\mathcal{N}, \mathcal{L})$ are finitely generated free left $D$-modules, then the map $\eta$ is natural with respect to morphisms $\mathcal{M} \rightarrow \mathcal{N}$, as one can check.

The properties we are interested in are as follows:  for all $i \in \mathbb{Z}$,
\begin{enumerate}
\item{} $\operatorname{End }\mathcal{L}_{i}= :D_{i}$ is a division ring with $k$ in its center,

\item{} $\Hom(\el_{i},\el_{i-1})=0=\operatorname{Ext}^{1}(\mathcal{L}_{i},\mathcal{L}_{i-1})$,

\item{} $\operatorname{Hom}(\mathcal{L}_{i}, \mathcal{L}_{i+1})$ is finite-dimensional as both a left $D_{i+1}$-module and a right $D_{i}$-module.  The left and right dimensions are denoted $l_{i}$ and $r_{i}$.

\item{} There is a short exact sequence, which we call the {\it $i$th Euler sequence},
$$
0 \rightarrow \mathcal{L}_{i} \overset{\eta_{\mathcal{L}_{i}}}{\rightarrow} {}^{*}\operatorname{Hom}(\mathcal{L}_{i}, \mathcal{L}_{i+1}) \otimes_{D_{i+1}} \mathcal{L}_{i+1} \rightarrow \mathcal{L}_{i+2} \rightarrow 0.
$$

\item{} The canonical one-to-one $k$-algebra map $\Phi_{i}:D_{i} \rightarrow D_{i+2}$ (defined in Lemma \ref{lemma.lstruct1} using (4)) is an isomorphism.

\item{} $\operatorname{Ext}^{1}(\mathcal{L}_{i}, \mathcal{L}_{j})=0$ for all $j \geq i$.

\item{} There is an equality $l_{i}=r_{i-1}$.

\item{} For all $\mathcal{M}$ in ${\sf C}$, $\operatorname{Hom}(\mathcal{L}_{i}, \mathcal{M})$ is a finite-dimensional right $D_{i}$-module.

\item{} The sequence $\underline{\mathcal{L}}$ is ample.
\end{enumerate}

Inspired by Seminaire Rudakov \cite{rudakov}, we make the following definition.
\begin{definition}
Let ${\sf C}$ be a $k$-linear abelian category.  A sequence of objects
$\underline{\mathcal{L}} := (\mathcal{L}_{i})_{i \in \mathbb{Z}}$ in ${\sf C}$ is called a {\it helix} if it satisfies properties (1)-(9).
\end{definition}
This terminology is also employed in \cite{newizuru} for a related notion.

We remark that, assuming (1) and (3) hold, property (4) is equivalent to the following:
there is a short exact sequence 
$$
0 \rightarrow \mathcal{L}_{i} \rightarrow \mathcal{L}_{i+1}^{ l_{i}} \rightarrow \mathcal{L}_{i+2} \rightarrow 0
$$
such that the $l_{i}$ morphisms defining the left arrow are left linearly independent.  This fact will be employed without comment in the sequel.

From the sequence $\underline{\mathcal{L}}$ we can form the $\mathbb{Z}$-algebra $H$ with $H_{ij}=\operatorname{Hom}(\el_{-j}, \el_{-i})$ and with multiplication induced by composition.

We prove the following result:

\begin{theorem} \label{theorem.supermain}
Let $\underline{\mathcal{L}}$ denote a sequence of objects in a $k$-linear abelian category ${\sf C}$.
\begin{enumerate}
\item{} If $\underline{\mathcal{L}}$ satisfies (1)-(5), then $H_{ii+1}$ is admissible and there is a $k$-linear $\mathbb{Z}$-algebra homomorphism
\begin{equation} \label{eqn.mainiso}
\mathbb{S}^{nc}(H_{01}) \rightarrow H
\end{equation}
which is an isomorphism in degrees zero and one.

\item{} If $\underline{\mathcal{L}}$ satisfies (1)-(6), then the homomorphism (\ref{eqn.mainiso}) is an epimorphism.

\item{} If $\underline{\mathcal{L}}$ satisfies (1)-(7), then the homomorphism (\ref{eqn.mainiso}) is an isomorphism, and $H_{ii+1}$ is 2-periodic and not of type $(1,1)$, $(1,2)$ or $(1,3)$.

\item{} If the category ${\sf C}$ has a helix $\underline{\mathcal{L}}$, then there is a $k$-linear equivalence ${\sf C} \equiv \mathbb{P}^{nc}(M)$ with $M$ a 2-periodic bimodule not of type $(1,1)$, $(1,2)$ or $(1,3)$ such that $\mathbb{S}^{nc}(M)$ is graded coherent.

\item{} If there is a $k$-linear equivalence ${\sf C} \equiv \mathbb{P}^{nc}(M)$ with $M$ a 2-periodic bimodule not of type $(1,1)$, $(1,2)$ or $(1,3)$ such that $\mathbb{S}^{nc}(M)$ is graded coherent, then ${\sf C}$ has a helix.
\end{enumerate}

\end{theorem}
By Theorem \ref{theorem.tsen} and Theorem \ref{theorem.supermain}, Kussin's noncommutative curves of genus zero \cite{kussin} have helices.  In addition, although the noncommutative curves studied in \cite{newchan} are noetherian, they are not noncommutative curves of genus zero in the sense of \cite{kussin} since they are not necessarily Ext-finite.  Nevertheless, Theorem \ref{theorem.supermain} implies they also have  helices.

We now describe some immediate consequences of some of these properties.

\begin{lemma}
Suppose $\underline{\mathcal{L}}$ satisfies (1)-(4).  Then for all $j<i$, $\Hom(\el_{i},\el_{j})=0$.
\end{lemma}

\begin{proof}
We proceed by induction on $i-j$.  The base case holds by (2). Now suppose $\operatorname{Hom}(\mathcal{L}_{i},\mathcal{L}_{j})=0$ for some $j < i$, and apply $\operatorname{Hom}(\mathcal{L}_{i},-)$ to the $j-1$th Euler sequence.  The induced sequence starts
$$
0 \rightarrow \operatorname{Hom}(\mathcal{L}_{i},\mathcal{L}_{j-1}) \rightarrow \operatorname{Hom}(\mathcal{L}_{i},\mathcal{L}_{j})^{ l_{j-1}}.
$$
Since the right term is zero so is $\operatorname{Hom}(\mathcal{L}_{i},\mathcal{L}_{j-1})$, as desired.
\end{proof}

\begin{lemma} \label{lemma.nonzero}
Suppose $\underline{\mathcal{L}}$ satisfies (1)-(4).  Then the numbers $l_{i}$ and $r_{i}$ are nonzero.
\end{lemma}

\begin{proof}
Suppose one of $l_{i}$ or $r_{i}$ was equal to zero.  Then $\operatorname{Hom}(\mathcal{L}_{i},\mathcal{L}_{i+1})=0$, so that the $i$th Euler sequence would imply that $\mathcal{L}_{i}=0$.  This contradicts the fact that $\operatorname{End}(\mathcal{L}_{i})=D_{i}$.
\end{proof}

\begin{lemma} \label{lemma.lstruct1}
Suppose $\underline{\mathcal{L}}$ satisfies (1), (3) and (4).

\begin{enumerate}
\item{} If $a \in D_{i}$, then there exists a map $f_{a} \in \operatorname{End}(\mathcal{L}_{i+1}^{ l_{i}})$ such that the diagram
$$
\begin{CD}
\mathcal{L}_{i} & \overset{h}{\longrightarrow} & \mathcal{L}_{i+1}^{ l_{i}}  \\
@V{a}VV  @VV{f_{a}}V \\
\mathcal{L}_{i} & \underset{h}{\longrightarrow} & \mathcal{L}_{i+1}^{ l_{i}}
\end{CD}
$$
whose horizontal arrows are from the $i$th Euler sequence, commutes.

\item{} The function $a \mapsto f_{a}$ is a $k$-algebra homomorphism.

\item{} There exists a $k$-algebra homomorphism
$$
\Phi_{i}:\operatorname{End}(\el_{i}) \longrightarrow \operatorname{End}(\el_{i+2})
$$
endowing $\operatorname{Hom}(\el_{i+1},\el_{i+2})$ with an $\operatorname{End}(\el_{i})-\operatorname{End}(\el_{i+1})$-bimodule structure.
\end{enumerate}
\end{lemma}

\begin{proof}
Properties (1) and (2) follow immediately from the fact that $\eta_{\mathcal{L}_{i}}$ is natural.

%By assumption, the components of $h=(h_{1},\ldots,h_{l_{i}})$ are a left basis for $\operatorname{Hom}(\mathcal{L}_{i},\mathcal{L}_{i+1})$.  Therefore, if $a \in D_{i}$, there exist $a_{lj}$ in $D_{i+1}$ such that
%\begin{equation} \label{eqn.lr}
%h_{l}a = \sum_{j} a_{lj}h_{j}.
%\end{equation}
%We let $f_{a} \in \operatorname{End}(\mathcal{L}_{i+1}^{ l_{i}})$ denote the morphism uniquely determined by the fact that if $g_{i}:\mathcal{L}_{i+1} \longrightarrow \mathcal{L}_{i+1}^{ l_{i}}$ denotes the $i$th inclusion, then $f_{a} g_{j} = (a_{1j},\ldots,a_{l_{i}j})$.  Uniqueness of $f_{a}$ follows from the fact that $\{h_{1},\ldots,h_{l_{i}}\}$ is a left basis for $\operatorname{Hom}(\mathcal{L}_{i},\mathcal{L}_{i+1})$.  Part (1) follows.

%The proof (2) is routine and omitted.

For (3), we define the $k$-algebra homomorphism $\Phi_{i}:D_{i} \longrightarrow D_{i+2}$ as follows:  given $a \in D_{i}$, part (1) implies that we get a map $f_{a} \in \operatorname{End}(\el_{i+1}^{l_{i}})$ such that the diagram
$$
\begin{CD}
\el_{i} & \overset{h}{\longrightarrow} & \el_{i+1}^{l_{i}} \\
@V{a}VV  @VV{f_{a}}V \\
\el_{i} & \underset{h}{\longrightarrow} & \el_{i+1}^{l_{i}}
\end{CD}
$$
commutes.  It follows that there is a unique $a' \in D_{i+2}$ making the diagram
\begin{equation} \label{eqn.ses}
\begin{CD}
\el_{i+1}^{l_{i}} & {\longrightarrow} & \el_{i+2}  \\
@VV{f_{a}}V @VV{a'}V \\
\el_{i+1}^{l_{i}} & {\longrightarrow} & \el_{i+2}
\end{CD}
\end{equation}
whose horizontal arrows are those in the $i$th Euler sequence, commute.  We define $\Phi_{i}(a) := a'$.  The proof that $\Phi_{i}$ is a $k$-algebra homomorphism is routine and omitted.
\end{proof}

\section{Sufficiency of the helical criterion} \label{section.abelian}
Throughout this section, we let ${\sf C}$ denote a $k$-linear abelian category, and we let $\underline{\mathcal{L}}$ be a sequence of objects in ${\sf C}$.  We say that $\underline{\mathcal{L}}$ is {\it linear} if it satisfies properties (1)-(5) from Section \ref{section.props}.  The purpose of this section is to study properties of linear sequences and use them to prove Theorem \ref{theorem.supermain}(1)-(4).  In particular, we show that the existence of a helix in ${\sf C}$ implies that ${\sf C}$ is a noncommutative projective line.

\subsection{Linear sequences}

\begin{prop} \label{prop.phi}
Suppose $\underline{\mathcal{L}}$ is linear.  If, for each $i$, we endow $\operatorname{Hom}(\el_{i+1},\el_{i+2})$ with the $D_{i} - D_{i+1}$-bimodule structure from Lemma \ref{lemma.lstruct1}(3), then there is an isomorphism of bimodules
$$
\Psi_{i}: {}^{*}\operatorname{Hom}(\el_{i},\el_{i+1}) \longrightarrow \operatorname{Hom}(\el_{i+1},\el_{i+2}).
$$
\end{prop}

\begin{proof}
The map $\Psi_{i}$ is constructed by applying $\operatorname{Hom}(\mathcal{L}_{i+1},-)$ to the $i$th Euler sequence
$$
0 \rightarrow \mathcal{L}_{i} \overset{\eta_{\mathcal{L}_{i}}}{\rightarrow} {}^{*}\operatorname{Hom}(\mathcal{L}_{i}, \mathcal{L}_{i+1}) \otimes_{D_{i+1}} \mathcal{L}_{i+1} \rightarrow \mathcal{L}_{i+2} \rightarrow 0
$$
and using a variant of \cite[Proposition B3.19]{az2}.  The map is an isomorphism by property (2).

\end{proof}

The following is an adaptation of \cite[Proposition 2.2]{dlab}.
\begin{lemma} \label{lemma.relations}
Suppose $\underline{\mathcal{L}}$ is linear, and suppose $\{h_{1}, \ldots, h_{l_{i}}\}$ is a left basis for  $\operatorname{Hom}(\el_{i},\el_{i+1})$.  Then, under the composition
\begin{eqnarray*}
{}^{*}\operatorname{Hom}(\el_{i},\el_{i+1}) \otimes \operatorname{Hom}(\el_{i},\el_{i+1}) & \overset{\Psi_{i} \otimes 1}{\longrightarrow} & \operatorname{Hom}(\el_{i+1},\el_{i+2}) \otimes \operatorname{Hom}(\el_{i},\el_{i+1}) \\
& \longrightarrow & \operatorname{Hom}(\el_{i},\el_{i+2})
\end{eqnarray*}
whose second arrow is induced by composition, the element $\sum_{j} {}^{*}h_{j} \otimes h_{j}$ goes to zero.
\end{lemma}

\begin{proof}
%As in the proof of Proposition \ref{prop.phi}, we let $g_{m}:\el_{i+1} \longrightarrow \el_{i+1}^{\oplus l_{i}}$ denote the $m$th inclusion.  By definition of $\Psi_{i}$, the element $\sum_{j} {}^{*}h_{j} \otimes h_{j}$ maps to $h' (\sum_{j} g_{j} h_{j})$.  

As one can check, the element $\sum_{j} {}^{*}h_{j} \otimes h_{j}$ maps to the composition of maps in the $i$th Euler sequence, hence maps to zero.
\end{proof}

\begin{prop} \label{prop.admissible}
Suppose $\underline{\mathcal{L}}$ is linear.  Then the $D_{i+1}-D_{i}$-bimodule $M: = \Hom(\el_{i}, \el_{i+1})$ is admissible.
\end{prop}

\begin{proof}
The fact that the left duals of $M$ are finite dimensional on either side follows immediately from Proposition \ref{prop.phi}.  On the other hand, if $N=\operatorname{Hom}(\el_{i-1},\el_{i})$ then Proposition \ref{prop.phi} implies that $M \cong {}^{*}N$.  Since $N$ is finite-dimensional on either side by linearity of $\underline{\mathcal{L}}$, the canonical map $N \rightarrow ({}^{*}N)^{*}$ is an isomorphism of bimodules, and so $M^{*} \cong N$.  Since $i$ is arbitrary, the result follows.
\end{proof}

In light of Lemma \ref{lemma.relations} and Proposition \ref{prop.admissible}, the proof of the next result is similar to the proof of \cite[Proposition 3.4]{tsen}.  We leave the details to the interested reader.

\begin{prop} \label{prop.ncsym1}
Suppose $\underline{\mathcal{L}}$ is linear, let $j \in \mathbb{Z}$, and let $M = \operatorname{Hom}(\el_{-(j+1)}, \el_{-j})$.  Then, for all $i \in \mathbb{Z}$, there is a canonical isomorphism of bimodules
$$
\Psi_{i}: M^{i*} \overset{\cong}{\rightarrow} \Hom(\el_{-j-(i+1)}, \el_{-j-i})
$$
where the bimodule structure on $\Hom(\el_{-j-(i+1)}, \el_{-j-i})$ is given by the appropriate composition of maps $\Phi_{l}$ defined in Lemma \ref{lemma.lstruct1}(3).  Furthermore,
$$
Q_{i} \subset \operatorname{ker }(\Psi_{i} \otimes \Psi_{i+1}).
$$
\end{prop}

\begin{cor} \label{cor.linear}
Suppose $\underline{\mathcal{L}}$ is linear, and let $i \in \mathbb{Z}$.  Then $H_{i,i+1}$ is admissible and there is a $k$-linear $\mathbb{Z}$-algebra homomorphism
\begin{equation} \label{eqn.mainhom}
\mathbb{S}^{nc}(H_{01}) \rightarrow H
\end{equation}
which is an isomorphism in degrees zero and one.
\end{cor}

\begin{proof}
By Proposition \ref{prop.admissible}, $H_{i,i+1}$ is admissible.  The rest of the result follows immediately from Proposition \ref{prop.ncsym1}.
\end{proof}

\subsection{Proof of Theorem \ref{theorem.supermain}(1)-(4)}
In this section we retain the notation from Section \ref{section.props}.  Our goal is to prove Theorem \ref{theorem.supermain}(1)-(4).  We begin by noting that Theorem \ref{theorem.supermain}(1) follows from Corollary \ref{cor.linear}.  Theorem \ref{theorem.supermain}(2) follows easily from the next result.

\begin{lemma}
Suppose $\underline{\mathcal{L}}$ is linear and satisfies property (6) from Section \ref{section.props}.  Then, for each $i \in \mathbb{Z}$ and $j>i+1$, the map
$$
\operatorname{Hom}(\mathcal{L}_{j-1}, \mathcal{L}_{j}) \otimes \cdots \otimes \operatorname{Hom}(\mathcal{L}_{i}, \mathcal{L}_{i+1}) \rightarrow \operatorname{Hom}(\mathcal{L}_{i}, \mathcal{L}_{j})
$$
induced by composition is surjective.
\end{lemma}

\begin{proof}
Let $j>i+1$ and apply $\operatorname{Hom}(\mathcal{L}_{i},-)$ to the $j-2$th Euler sequence.  Since $\operatorname{Ext}^{1}(\mathcal{L}_{i},\mathcal{L}_{j-2})=0$ by property (6) of $\underline{\mathcal{L}}$, the induced map $\operatorname{Hom}(\mathcal{L}_{i}, \mathcal{L}_{j-1}^{ l_{j-2}}) \rightarrow \operatorname{Hom}(\mathcal{L}_{i}, \mathcal{L}_{j})$ is surjective.  Therefore, the map
\begin{equation} \label{eqn.indmap}
\Hom(\el_{j-1}, \el_{j}) \otimes \Hom(\el_{i}, \el_{j-1}) \rightarrow \Hom(\el_{i}, \el_{j})
\end{equation}
induced by composition is surjective.

We now proceed to prove the lemma by induction on $j$.  If $j=i+2$, surjectivity of (\ref{eqn.indmap}) implies the base case.  If $j>i+2$, surjectivity of (\ref{eqn.indmap}) together with the induction hypothesis implies the result.
\end{proof}

\begin{proposition} \label{prop.forbid}
If $\underline{\mathcal{L}}$ is linear and satisfies (6) and (7), and $j \in \mathbb{Z}$, then
\begin{enumerate}
\item{} $H_{j,j+1}$ is 2-periodic, and

\item{} $H_{j,j+1}$ is not of type $(1,1)$, $(1,2)$ or $(1,3)$.
\end{enumerate}
\end{proposition}

\begin{proof}
To prove the first assertion, we note that $H_{j,j+1}$ is admissible by Corollary \ref{cor.linear}.  Furthermore
\begin{eqnarray*}
\operatorname{ldim}H_{j,j+1}^{i*} & = & \operatorname{ldim}\operatorname{Hom}(\el_{-j-(i+1)},\el_{-j-i}) \\
& = & l_{-j-(i+1)} \\
& = & r_{-j-(i+2)} \\
& = & \operatorname{rdim }\operatorname{Hom}(\el_{-j-(i+2)}, \el_{-j-(i+1)}) \\
& = & \operatorname{rdim }H_{j,j+1}^{i+1*}
\end{eqnarray*}
where the first and last equality follow from Proposition \ref{prop.ncsym1}, while the third equality follows from property (7) of $\underline{\mathcal{L}}$.  Therefore, $H_{j,j+1}$ is 2-periodic.

We next show that if $\underline{\mathcal{L}}$ satisfies properties (1)-(7), then $H_{01}$ is not of type $(1,3)$.  By 2-periodicity of $H_{01}$, it will follow that $H_{j,j+1}$ is not of type $(1,3)$.  The rest of the proof is similar (but easier) and omitted.

Suppose $H_{01}$ is of type $(1,3)$.  We first assume $l_{-1}=1$ while $r_{-1}=3$.  Then by 2-periodicity of $H_{01}$ and Proposition \ref{prop.ncsym1}, if $i$ has odd parity then $l_{i}=1$ and $r_{i}=3$, while if $i$ has even parity, then $l_{i}=3$ and $r_{i}=1$.

Suppose, first, that $i$ has odd parity.  Then the $-i$th Euler sequence has the form
$$
0 \rightarrow \mathcal{L}_{-i} \rightarrow \mathcal{L}_{-i+1} \rightarrow \mathcal{L}_{-i+2} \rightarrow 0.
$$
Applying $\operatorname{Hom}(\mathcal{L}_{-i},-)$, we deduce, by property (6) of $\underline{\mathcal{L}}$, that
\begin{equation} \label{eqn.rtdim}
\operatorname{rdim }\operatorname{Hom}(\mathcal{L}_{-i},\mathcal{L}_{-i+2})=2.
\end{equation}
Next, applying $\operatorname{Hom}(\mathcal{L}_{-i},-)$ to the $-i+1$st Euler sequence
$$
0 \rightarrow \mathcal{L}_{-i+1} \rightarrow \mathcal{L}_{-i+2}^{3} \rightarrow \mathcal{L}_{-i+3} \rightarrow 0
$$
and using (\ref{eqn.rtdim}) allows us to deduce, by property (6) of $\underline{\mathcal{L}}$, that
$$
\operatorname{rdim }\operatorname{Hom}(\mathcal{L}_{-i},\mathcal{L}_{-i+3})=3.
$$
Similarly, $\operatorname{rdim }\operatorname{Hom}(\mathcal{L}_{-i},\mathcal{L}_{-i+4})=1$ and thus $\operatorname{rdim }\operatorname{Hom}(\mathcal{L}_{-i},\mathcal{L}_{-i+5})=0$.  But now applying $\operatorname{Hom}(\mathcal{L}_{-i},-)$ to
$$
0 \rightarrow \mathcal{L}_{-i+4} \rightarrow \mathcal{L}_{-i+5} \rightarrow \mathcal{L}_{-i+6} \rightarrow 0.
$$
allows us to deduce that $\operatorname{Hom}(\mathcal{L}_{-i},\mathcal{L}_{-i+4})$ is zero, a contradiction.

One obtains a similar contradiction in the case that $i$ has even parity.
\end{proof}

In light of Proposition \ref{prop.forbid}, Theorem \ref{theorem.supermain}(3) and Theorem \ref{theorem.supermain}(4) will follow from the next

\begin{prop} \label{prop.isoo}
If $\underline{\mathcal{L}}$ is linear and satisfies (6) and (7), then the homomorphism
$$
\mathbb{S}^{nc}(H_{01}) \rightarrow H
$$
from Corollary \ref{cor.linear} is an isomorphism.

If $\underline{\mathcal{L}}$ is a helix then the isomorphism above induces an equivalence
$$
\mathbb{P}^{nc}(H_{01}) \rightarrow {\sf C}.
$$
\end{prop}

\begin{proof}
To prove the first result, it suffices, by Theorem \ref{theorem.supermain}(2), to show that for all $i,j \in \mathbb{Z}$, the right dimension of $H_{ij}$ equals that of $A_{ij} := \mathbb{S}^{nc}(H_{01})$.  By definition of the noncommutative symmetric algebra, and by properties (2) and (4) of $\underline{\mathcal{L}}$, $H_{ij}$ and $A_{ij}$ have right-dimension zero if $j<i$. Furthermore, by Corollary \ref{cor.linear}, the result holds if $j-i$ equals $0$ or $1$.  Now we prove the result by induction on the difference of indices.  By properties (4) and (6) of $\underline{\mathcal{L}}$, we have an exact sequence of right $D_{-j}$-modules
$$
0 \rightarrow H_{ij} \rightarrow H_{i-1, j}^{\oplus l_{-i}} \rightarrow H_{i-2, j} \rightarrow 0.
$$
On the other hand, by Corollary \ref{cor.euler}, we have an exact sequence of right-$A_{jj}\cong D_{-j}$-modules
$$
0 \rightarrow A_{ij} \rightarrow A_{i-1, j}^{\oplus \mbox{rdim }A_{i-2, i-1}} \rightarrow A_{i-2, j} \rightarrow 0
$$
since $j \geq i+2$.  Now, $l_{-i}=r_{-i-1}$ by (7), while by Proposition \ref{prop.phi}, $l_{-i}=r_{-i+1}$, so that $r_{-i-1}$ is the right-dimension of $\operatorname{Hom}(\el_{-i+1}, \el_{-i+2})=H_{i-2 i-1}$.  Finally, by Proposition \ref{prop.ncsym1}, $H_{i-2, i-1} \cong A_{i-2, i-1}$. The first result now follows by induction.

Now we prove the second result.  Suppose $\underline{\mathcal{L}}$ is a helix.  Since (8) and (9) hold, the indicated equivalence follows from the first part and Theorem \ref{theorem.polish}.
\end{proof}

\subsection{The finite over $k$ case}
We say $\underline{\mathcal{L}}$ is {\it finite over $k$} if it satisfies (1) and each $D_{i}$ is finite-dimensional over $k$.

Theorem \ref{theorem.ft} follows from the next
\begin{lemma} \label{lemma.ft}
Suppose $\underline{\mathcal{L}}$ is finite over $k$ and $\operatorname{dim }_{k}D_{i}=\operatorname{dim }_{k}D_{i+2}$ for all $i$.  Suppose further that $\underline{\mathcal{L}}$ satisfies properties (2), (4), and (8).  Then $\underline{\mathcal{L}}$ satisfies properties (3), (5), and (7).
\end{lemma}

\begin{proof}
By property (8), $\operatorname{Hom}(\el_{i}, \el_{i+1})$ is finite-dimensional on the right.  Hence, it is finite-dimensional over $k$ on the left and therefore finite-dimensional over $D_{i+1}$ on the left so that (3) holds.

Property (5) holds since, by Lemma \ref{lemma.lstruct1}(3), $\Phi_{i}$ is a $k$-algebra injection and, by hypothesis, $\operatorname{dim }_{k}D_{i}=\operatorname{dim }_{k}D_{i+2}$.

Finally, we check (7).  Let $m_{i}$ denotes the dimension of $D_{i}$ over $k$.  Then
\begin{eqnarray*}
r_{i}m_{i} & = & l_{i}m_{i+1} \\
& = & r_{i+1}m_{i+1} \\
& = & l_{i+1}m_{i+2},
\end{eqnarray*}
where the second equality follows from Proposition \ref{prop.phi}.  The result follows.
\end{proof}
The next result follows immediately.

\begin{cor} \label{cor.central}
Suppose $\underline{\mathcal{L}}$ is such that $\operatorname{End}(\mathcal{L}_{i})=k$ for all $i$, and satisfies (2), (4) and (8). Then $\underline{\mathcal{L}}$ satisfies (3), (5) and (7).
\end{cor}

\section{Internal Hom functors and their derived functors} \label{section.shom}
To complete the proof of Theorem \ref{theorem.supermain}, it remains to prove Theorem \ref{theorem.supermain}(5).  To this end, we show that noncommutative projective lines are homologically well-behaved.  More specifically, in Section \ref{section.gor}, we prove that if $D_{0}$ and $D_{1}$ are division rings over $k$, $M$ is a $k$-central $D_{0}-D_{1}$-bimodule such that $M$ is 2-periodic and not of type $(1,1)$, $(1,2)$ or $(1,3)$, then $\mathbb{S}^{nc}(M)$ is Gorenstein (Corollary \ref{cor.gor0}).  This will allow us to prove homological results about $\mathbb{P}^{nc}(M)$, Corollary \ref{cor.gorenstein}, Theorem \ref{thm.vanish}, and Lemma \ref{lemma.sofinite} that we will need to prove Theorem \ref{theorem.supermain}(5).

In this section we will describe many of the homological preliminaries we will need in Section \ref{section.gor}: in Section \ref{section.inthom} we recall the definition and basic properties of the internal Hom functor introduced in \cite{duality}.  Since we work here over an affine base, the constructions in \cite{duality} simplify considerably.  For this reason, and for the convenience of the reader, we will develop the material from first principals.  Next, in Section~\ref{section.intext}, we study the right-derived functors of the internal Hom functor.

Throughout this section, if $R$ is a ring, ${\sf Mod }R$ will denote the category of right $R$-modules.  Furthermore, unless otherwise stated, $A$ will denote a $\mathbb{Z}$-algebra such that, for all $i$, $A_{ii}$ is a division ring with $k$ in its center.  To define the internal Hom functor, we will need the following

\begin{definition}
Let ${\sf Bimod }A-A$ denote the category defined as follows:
\begin{itemize}

\item{}
an object of ${\sf Bimod }A-A$ is a triple
$$
(C=\{C_{ij}\}_{i,j \in \mathbb{Z}}, \{\mu_{ijk}\}_{i,j,k \in \mathbb{Z}}, \{\psi_{ijk}\}_{i,j,k \in \mathbb{Z}})
$$
where $C_{ij}$ is an $A_{ii}-A_{jj}$-bimodule and $\mu_{ijk}:C_{ij} \otimes A_{jk} \rightarrow C_{ik}$ and $\psi_{ijk}: A_{ij} \otimes C_{jk} \rightarrow C_{ik}$ are morphisms of $A_{ii}-A_{kk}$-bimodules making $C$ both a graded right and left $A$-module such that the actions are compatible in the usual sense.

\item{}
A morphism $\phi:  C \rightarrow D$ between objects in ${\sf Bimod }A-A$ is a  collection $\{\phi_{ij}\}_{i,j \in \mathbb{Z}}$, where $\phi_{ij}:C_{ij} \rightarrow D_{ij}$ are morphisms of $A_{ii}-A_{jj}$-bimodules which respect the $A-A$-bimodule structure on $C$ and $D$.
\end{itemize}
\end{definition}
We omit the routine verification that ${\sf Bimod }A-A$ is abelian.

\subsection{Internal Hom} \label{section.inthom}
We begin by defining the internal Hom functor.  For $C$ an object in ${\sf Bimod }A-A$ and $M$ an object in ${\sf Gr }A$,
\begin{itemize}
\item{} we let
$$
\sHom(e_{i}C,M)
$$
denote the right $A_{ii}$-module with underlying set  $\operatorname{Hom}_{{\sf Gr}A}(e_{i}C,M)$ and with $A_{ii}$-action induced by the left action of $A_{ii}$ on $e_{i}C$, and

\item{} we let
$$
{\usHom}(C,M)
$$
denote the object in ${\sf Gr }A$ with $i$th component $\sHom(e_{i}C,M)$ and with multiplication induced by left-multiplication of $A$ on $C$.
\end{itemize}

\begin{lemma} \label{lemma.basic}
The assignment
$$
\usHom(-,-):({{\sf Bimod }A-A})^{op} \times {\sf Gr }A \rightarrow {\sf Gr}A
$$
is a bifunctor in the sense of \cite[Chapter 20.8]{lang}.  Furthermore,
\begin{enumerate}
\item{} for each $M$ in ${\sf Gr }A$, there is a natural isomorphism $\usHom(A,M) \cong M$, and

\item{} $\usHom(-,-)$ is left exact in each variable, and for each injective $I$ in ${\sf Gr }A$, $\usHom(-,I)$ is exact.
\end{enumerate}
\end{lemma}

\begin{proof}
The fact that $\usHom(-,-)$ is a bifunctor follows easily from the fact that $\operatorname{Hom}_{{\sf Gr }A}(-,-)$ is a bifunctor.  The proof of the first item is routine while the second follows from properties of $\operatorname{Hom}_{{\sf Gr }A}(-,-)$.
\end{proof}

If $C$ is an object of ${\sf Bimod }A-A$, then we define two associated functors
$$
-\otimes_{A_{ii}}e_{i}C: {\sf Mod }A_{ii} \rightarrow {\sf Gr }A,
$$
and
$$
\sHom(e_{i}C,-):{\sf Gr }A \rightarrow {\sf Mod }A_{ii}.
$$
As one can check, the pair $(-\otimes_{A_{ii}}e_{i}C, \sHom(e_{i}C,-))$ has a canonical adjoint structure.  This is exploited in the proof of the following

\begin{theorem} \label{theorem.local}
Suppose $F$ is an $A_{jj}-A_{ii}$-bimodule of finite dimension on either side, and let $F \otimes_{A_{ii}} e_{i}C$ denote the object of ${\sf Bimod }A-A$ such that
$$
(F \otimes_{A_{ii}} e_{i}C)_{lm} = \begin{cases} F \otimes_{A_{ii}}C_{im} & \mbox{if $l=j$} \\ 0 & \mbox{otherwise,} \end{cases}
$$
endowed with the obvious bimodule structure.  Let $M$ be an object of ${\sf Gr} A$.  Then

\begin{enumerate}
\item{}
there is a natural isomorphism of ${\sf Mod }A_{jj}$-valued functors
$$
\sHom(F \otimes e_{i}C,-) \cong \sHom(e_{i}C,-) \otimes_{A_{ii}} F^{*},
$$
and

\item{}
the functor $\sHom(F \otimes e_{i}A,-):{\sf Gr }A \rightarrow {\sf Mod }A_{jj}$ is exact.
\end{enumerate}
\end{theorem}

\begin{proof}
By adjointness, there is a canonical isomorphism
$$
\operatorname{Hom}_{{\sf Gr }A}(F \otimes e_{i}C,M) \rightarrow \operatorname{Hom}_{{\sf Mod }A_{ii}}(F, \sHom(e_{i}C,M))
$$
where $F$ is considered as an $A_{ii}$-module.  Furthermore, as one can check, this isomorphism is compatible with right $A_{jj}$-module structure.  Finally, by the Eilenberg-Watts theorem, there is a natural isomorphism
$$
\operatorname{Hom}_{{\sf Mod }A_{ii}}(F,\sHom(e_{i}C,M)) \rightarrow \sHom(e_{i}C,M) \otimes_{A_{ii}}F^{*}
$$
of right $A_{jj}$-modules, completing the proof of the first assertion.

The proof of the second assertion follows from the first, in light of the fact that the functor $-\otimes_{A_{ii}} F^{*}: {\sf Mod }A_{ii} \rightarrow {\sf Mod }A_{jj}$ is exact (see, for example, \cite[Section 2.1]{witt}), and that $\sHom(e_{i}A,-) \cong (-)_{i}$.

\end{proof}

\subsection{Internal Ext} \label{section.intext}
Let $C \in {\sf Bimod }A-A$.  In this section we study the right derived functors of $\usHom(C,-)$ and $\sHom(e_{j}C,-)$.  The fact that $\usHom(C,-)$ and $\sHom(e_{j}C,-)$ have right derived functors follows from Lemma \ref{lemma.basic}(2).  We denote them by $\usExt^{i}(C,-)$ and $\sExt^{i}(e_{j}C,-)$.  We note that since taking the $j$th degree part of an object of ${\sf Gr }A$ is an exact functor from ${\sf Gr }A$ to ${\sf Mod }A_{jj}$, we have
$$
(\usExt^{i}(C,M))_{j} \cong \sExt^{i}(e_{j}C,M).
$$

In order to state the next lemma, we need to introduce some terminology.  We say an object $C$ of ${\sf Bimod }A-A$ is {\it left-bounded by degree $l$} if for each $i$, $C_{ij} \neq 0$ implies $j \geq i+l$.  For $n$ a nonnegative integer, we let $A_{\geq n}$ denote the subobject of $A$ in ${\sf Bimod }A-A$ given by $\oplus_{j-i \geq n} A_{ij}$.

\begin{lemma} \label{lemma.basic2}
For $C \in {\sf Bimod }A-A$ and $M \in {\sf Gr }A$,
\begin{enumerate}
\item{} the sequence $\usExt^{i}(-,M)$ forms a $\delta$-functor,

\item{} if $M$ is right-bounded by degree $r$ and $C$ is left-bounded by degree $l$, then $\usExt^{i}(C,M)$ is right-bounded by degree $r-l$, and

\item{} if $M$ is right-bounded then for all $j \geq 1$, the right bound of $\usExt^{j}(A/A_{\geq n},M)$ tends to $-\infty$ as $n \rightarrow \infty$.
\end{enumerate}
\end{lemma}

\begin{proof}
The first result follows directly from \cite[Proposition 8.4, p. 810]{lang} in light of Lemma \ref{lemma.basic}(2).

To prove the second result, we adapt the proof of \cite[Proposition 3.1(2)]{az} to our context.  To this end, we claim that $\sExt^{i}(e_{j}C,M)=0$ if and only if $\operatorname{Ext}_{{\sf Gr }A}^{i}(e_{j}C,M)=0$.  This follows from the fact that if $i_{*}:{\sf Mod }A_{jj} \rightarrow {\sf Mod }k$ is the restriction of scalars functor, then $i_{*} \sHom(e_{j}C,-) \cong \operatorname{Hom}_{{\sf Gr }A}(e_{j}C,-)$. The claim implies that the argument in the proof of \cite[Proposition 3.1(2)]{az} may be used to prove the second result.

The third result is the $\mathbb{Z}$-algebra version of \cite[Proposition 3.1(5)]{az}, and we recount the proof for the convenience of the reader.  By part (1), the functor $\usHom(-,M)$ applied to the short exact sequence in ${\sf Bimod }A-A$
\begin{equation} \label{equation.sesa}
0 \rightarrow A_{\geq n} \rightarrow A \rightarrow A/A_{\geq n} \rightarrow 0
\end{equation}
induces a long-exact sequence.  Since, by Lemma \ref{lemma.basic}(1), $\usHom(A,-)$ is exact, it thus suffices to show that, for $i \geq 0$, $\usExt^{i}(A_{\geq n},M)$ is right-bounded by $r-n$, where $r$ is the right-bound of $M$.  This is exactly part (2) so the result follows.
\end{proof}

\begin{lemma} \label{lemma.extcommute}
Suppose $F$ is an $A_{jj}-A_{ii}$-bimodule of finite dimension on either side, and let $F \otimes_{A_{ii}} e_{i}C$ denote the object of ${\sf Bimod }A-A$ defined in Theorem \ref{theorem.local}.  If $M$ be an object of ${\sf Gr} A$, then there is a natural isomorphism of ${\sf Mod }A_{jj}$-valued functors
$$
\sExt^{j}(F \otimes e_{i}C,-) \cong \sExt^{j}(e_{i}C,-) \otimes F^{*}.
$$
\end{lemma}

\begin{proof}
The $j=0$ case follows from Theorem \ref{theorem.local}(1).

By Lemma \ref{lemma.basic2}(1), and the fact that the composition of a $\delta$-functor with an exact functor is a $\delta$-functor, the sequence $\sExt^{j}(-,M)$ is a $\delta$-functor.  Furthermore, for $M$ injective and $j>0$, $\sExt^{j}(-,M)=0$ by Lemma \ref{lemma.basic}(2).  The result now follows from \cite[Theorem 1.3A, p. 206]{hartshorne}.
\end{proof}
Suppose $\lambda, \rho \in \mathbb{Z}$ with $\lambda \leq \rho$ and let $M \in {\sf Gr }A$.  We write $M \subset [\lambda, \rho]$ if $M_{i}$ nonzero implies that $\lambda \leq i \leq \rho$.  Similarly, suppose $l,r \in \mathbb{Z}$ with $l \leq r$ and let $C \in {\sf Bimod }A$.  We write $C \subset [l,r]$ if $C_{ij}$ nonzero implies $l \leq j-i \leq r$.  We say $C$ is {\it concentrated in degree $m$} if $C \subset [m,m]$.  We let $A_{0}$ denote the quotient $A/A_{\geq 1}$ in ${\sf Bimod }A-A$.

\begin{lemma} \label{lemma.extbound}
Let $i$ be a nonnegative integer, let $C \in {\sf Bimod }A$ be such that $C \subset [l,r]$, and let $M \in {\sf Gr }A$ be such that
$$
\usExt^{j}(A_{0},M) \subset [\lambda, \rho]
$$
for all $j \geq i$.  Then, for $j \geq i$,
$$
\usExt^{j}(C,M) \subset [\lambda-r,\rho-l].
$$
\end{lemma}

\begin{proof}
First, assume $C$ is concentrated in degree $m$.  By Lemma \ref{lemma.extcommute},
$$
\usExt^{i}(C,M)_{n} \cong \usExt^{i}(A_{0},M)_{n+m} \otimes C^{*}_{n, n+m}.
$$
Thus, since $\usExt^{j}(A_{0},M) \subset [\lambda, \rho]$ for $j \geq i$, we have $\usExt^{j}(C,M) \subset [\lambda-m,\rho-m]$ for $j \geq i$.

Now suppose $C \subset [l,r]$ and define a subobject, $C'$ of $C$ by letting $e_{n}C'= C_{n, n+r}$. Then we have an exact sequence
\begin{equation} \label{eqn.holast}
0 \rightarrow C' \rightarrow C \rightarrow C/C' \rightarrow 0.
\end{equation}
Furthermore, by construction, $C/C' \subset [l, r-1]$.  We now prove the result by induction on $r-l$, the case $r=l$ being proven above.

Suppose the result holds when $r-l<m$ and let $C \in {\sf Bimod }A$ be such that $C \subset [l,r]$, with $r-l=m$.  The long exact sequence for $\usHom(-,M)$ applied to (\ref{eqn.holast}) contains
$$
\usExt^{j}(C/C',M) \rightarrow \usExt^{j}(C,M) \rightarrow \usExt^{j}(C',M).
$$
If $j \geq i$, then the left module is contained in $[\lambda-r+1, \rho-l]$, whereas the right module is in $[\lambda-r, \rho-r]$, so that the assertion follows.
\end{proof}

For the remainder of this section, we assume $A$ is a connected $\mathbb{Z}$-algebra, finitely generated in degree one, so that ${\sf Proj }A$, as well as the functors $\pi, \tau,$ and $\omega$, are defined.

\begin{lemma} \label{lemma.torsfunct}
There is an isomorphism of functors $\tau(-) \cong \underset{n \to \infty}{\lim} \usHom(A/A_{\geq n}, -)$.
\end{lemma}

\begin{proof}
Let $M \in {\sf Gr }A$.  By Lemma \ref{lemma.basic}(1) the canonical map $A \rightarrow A/A_{\geq n}$ induces an inclusion
$$
\dlim \usHom(A/A_{\geq n},M) \rightarrow M
$$
natural in $M$.  It is a routine verification to prove that the image of this map is exactly $\tau(M)$ and the result follows easily from this.
\end{proof}

\begin{theorem} \label{theorem.quotients}
Suppose that for all $n \geq 0$ and all $j \in \mathbb{Z}$, $e_{j}A/e_{j}A_{\geq n}$ has a finite resolution by objects $C \in {\sf Bimod }A-A$ such that $\sHom(e_{j}C,-)$ is exact.  Let $N$ be a graded $A$-module.  For $j \geq 1$, we have
$$
\mbox{R}^{j}\omega(\pi(N)) \cong \operatorname{\dlim} \usExt^{j}(A_{\geq n},N).
$$
Furthermore, for $i \geq 1$, the right-derived functors of $\tau$ and $\omega$ satisfy
$$
\mbox{R}^{i+1}\tau(N) \cong \mbox{R}^{i}\omega (\pi(N)),
$$
and, for each $N \in {\sf Gr }A$, there is an exact sequence
$$
0 \rightarrow \tau N \rightarrow N \rightarrow \omega \pi N \rightarrow \mbox{R}^{1}\tau N \rightarrow 0,
$$
whose central arrow is natural.
\end{theorem}

\begin{proof}
We adapt the proof of \cite[Lemma 4.1.5 and Lemma 4.1.6]{bv} to our context.  We first claim, as in \cite[Lemma 4.1.3]{bv}, that for $T \in {\sf Gr }A$ a torsion module, we have $\mbox{R}^{i}\tau(T)=0$ for $i>0$.  To prove the claim, it suffices, by Lemma \ref{lemma.torsfunct}, to prove that
$$
\dlim \sExt^{i}(e_{j}A/e_{j}A_{\geq n},T)=0
$$
for $i>0$ and all $j$.  By \cite[Proposition 8.2, p. 809]{lang}, $\sExt^{i}(e_{j}A/e_{j}A_{\geq n},-)$ can be computed as the $n$th cohomology of $\sHom(-,-)$ applied to the finite resolution of $e_{j}A/e_{j}A_{\geq n}$ which exists by hypothesis.  Therefore, the functor
$$
\dlim \sExt^{i}(e_{j}A/e_{j}A_{\geq n},-)
$$
commutes with direct limits, so that we may assume, without loss of generality, that $T$ is right-bounded by degree $r$.  The claim now follows from Lemma \ref{lemma.basic2}(3).

We next claim that $\omega \pi N \cong \dlim \usHom(A_{\geq n},N)$.  To prove this, we note that by the first claim, and by the fact that $A$ is a connected $\mathbb{Z}$-algebra, finitely generated in degree one, we can copy the proof of \cite[Lemma 4.1.4]{bv}.

Now we prove the theorem.  The first statement follows from the second claim, together with the fact, proven in \cite[Lemma 4.1.6]{bv}, that $R^{i}(\omega \pi) \cong R^{i}\omega \circ \pi$.  The second statement comes the long exact sequence constructed by applying
$$
\dlim \usHom(-,N)
$$
to the short exact sequence (\ref{equation.sesa}).
\end{proof}
If $D_{0}$ and $D_{1}$ are division rings over $k$, and $A = \mathbb{S}^{nc}(M)$ where $M$ is a 2-periodic $D_{0}$-$D_{1}$-bimodule over $k$ not of type $(1,1)$, $(1,2)$ or $(1,3)$, then $A$ is a connected $\mathbb{Z}$-algebra, finitely generated in degree one.  Furthermore, it follows from Corollary \ref{cor.euler} and Theorem \ref{theorem.local}(2) that for all $n \geq 0$ and all $j \in \mathbb{Z}$, $e_{j}A/e_{j}A_{\geq n}$ has a finite resolution by objects $C \in {\sf Bimod }A-A$ such that $\sHom(e_{j}C,-)$ is exact.  Therefore, the hypotheses of Theorem \ref{theorem.quotients} hold for $A$.

\section{Relative local cohomology of noncommutative symmetric algebras} \label{section.gor}
In this section we show that noncommutative symmetric algebras are Gorenstein (Theorem \ref{thm.gorenstein}), and apply this result to the computation of the right derived functors of the torsion functor $\tau$ (Corollary \ref{cor.gor0}).  This allows us to compute certain cohomology groups over $\mathbb{P}^{nc}(M)$, and these computations are used to complete the proof of Theorem \ref{theorem.supermain} in Section \ref{section.proof2}.

We assume throughout this section that $D_{0}$ and $D_{1}$ are division rings over $k$, $M$ is a 2-periodic $D_{0}$-$D_{1}$ $k$-central bimodule not of type $(1,1)$, $(1,2)$ or $(1,3)$, and $A = \mathbb{S}^{nc}(M)$.

\subsection{$\mathbb{S}^{nc}(M)$ is Gorenstein}
By the remarks preceding Theorem \ref{theorem.quotients}, the derived functors of $\sHom(A/A_{\geq 1},-)$ and $\usHom(A/A_{\geq 1},-)$ may be computed using the Euler sequence from Corollary \ref{cor.euler}.  This fact will be utilized in the proof of the next

\begin{theorem} \label{thm.gorenstein}
Let $i \geq 0$ and let $l$ and $j$ be integers.  Then
$$
\usExt^{i}(A/A_{\geq 1}, e_{l} A)=0 \mbox{ for $i \neq 2$}
$$
and
$$
\usExt^{2}(A/A_{\geq 1}, e_{l} A)_{j}  \cong
\begin{cases}
A_{l-2,l-2} & \text{if $j = l -2$}, \\
0& \text{otherwise}.
\end{cases}
$$
\end{theorem}

\begin{proof}
By the comment preceding the theorem, we may compute
$$
\usExt^{i}(A/A_{\geq 1}, e_{l}A)_{m}
$$
by taking the cohomology of the sequence
\begin{equation} \label{eqn.gorcomp}
\sHom(e_{m}A, e_{l}A) \overset{d_{0}}{\rightarrow} \sHom(A_{m, m+1} \otimes e_{m+1}A, e_{l}A) \overset{d_{1}}{\rightarrow} \sHom(Q_{m} \otimes e_{m+2}A,e_{l}A)
\end{equation}
coming from the application of $\sHom(-,e_{l}A)$ to the truncation of the exact sequence
$$
Q_{m}\otimes e_{m+2}A \rightarrow A_{m,m+1} \otimes e_{m+1}A \rightarrow e_{m}A
$$
from Corollary \ref{cor.euler}.  In particular, we have
$$
\usExt^{0}(A/A_{\geq 1},e_{l}A)_{m}=\operatorname{ker }d_{0},
$$
$$
\usExt^{1}(A/A_{\geq 1},e_{l}A)_{m}=\frac{\operatorname{ker }d_{1}}{\operatorname{im }d_{0}},
$$
$$
\usExt^{2}(A/A_{\geq 1},e_{l}A)_{m}=\frac{\sHom(Q_{m} \otimes e_{m+2}A,e_{l}A)}{\operatorname{im}d_{1}}
$$
and $\usExt^{i}(A/A_{\geq 1},e_{l}A)_{m}=0$ for $i>2$.

If $m<l-2$ it follows from Theorem \ref{theorem.local} that all terms in (\ref{eqn.gorcomp}) are zero so that the indicated groups are zero.  Similarly, if $m=l-2$, the first two terms of (\ref{eqn.gorcomp}) vanish.  Therefore, $\usExt^{0}(A/A_{\geq 1},e_{l}A)_{l-2}=\usExt^{1}(A/A_{\geq 1},e_{l}A)_{l-2}=0$ and
$$
\usExt^{2}(A/A_{\geq 1},e_{l}A)_{l-2} \cong A_{ll} = A_{l-2,l-2}
$$
by Theorem \ref{theorem.local}.  Therefore, to establish the theorem, we must prove that the sequence (\ref{eqn.gorcomp}) is exact for $m>l-2$.

We first show $\operatorname{ker }d_{0}=0$.  If $m=l-1$, this holds since the left-most term of (\ref{eqn.gorcomp}) is zero.  Next, suppose $f:e_{m}A \rightarrow e_{l}A$, and let
$$
\mu: A_{m,m+1} \otimes e_{m+1}A \rightarrow e_{m}A
$$
denote multiplication.  If $m=l$, then either $f=0$ or $f$ is an isomorphism.  In the latter case, since multiplication $\mu$ is nonzero, $d_{0}(f) \neq 0$.  Therefore, the result holds in this case.  Thus, we suppose $m>l$.  If $d_{0}(f)=0$, then since $\mu$ is surjective in degree $\geq m+1$, it suffices to show that if $f$ is left-multiplication by $x \in A_{lm}$, then $xy =0$ for all $y \in A_{m,m+1}$ implies that $x = 0$.  This follows from Lemma \ref{lemma.euler}.

Next, we prove that if $m>l-2$, then $\operatorname{ker }d_{1}=\operatorname{im }d_{0}$.  It suffices to prove that $\operatorname{ker }d_{1} \subset \operatorname{im }d_{0}$.  Suppose $g:A_{m,m+1} \otimes e_{m+1}A \rightarrow e_{l}A$ is in $\operatorname{ker }d_{1}$.  Then $g$ has a factorization
$$
A_{m,m+1} \otimes e_{m+1}A \overset{\mu}{\rightarrow} (e_{m}A)_{\geq m+1} \overset{f}{\rightarrow} e_{l}A.
$$
Thus, to complete the proof of the result in this case, we must show that $f$ extends to a right $A$-module map $\tilde{f}:e_{m}A \rightarrow e_{l}A$.  Suppose $x_{1}, \ldots, x_{n}$ is a right basis for $A_{m, m+1}$ with associated right duals $x_{1}^{*}, \ldots, x_{n}^{*} \in A_{m+1,m+2}$.  Then
\begin{eqnarray*}
\sum_{p}f(x_{p})x_{p}^{*} & = & f(\sum_{p}x_{p}x_{p}^{*}) \\
& = & 0.
\end{eqnarray*}
It thus follows from Theorem \ref{theorem.euler} that $\sum_{p}f(x_{p})\otimes x_{p}^{*} \in A_{l ,m+1} \otimes A_{m+1,m+2}$ is an element of the image of $A_{lm} \otimes Q_{m}$ under the left map of (\ref{eqn.seqq}).   Therefore, there exists a $y \in A_{lm}$ such that, for all $p$, $f(x_{p})=y x_{p}$.  We let $\tilde{f}$ be defined by $\tilde{f}(1)=y$, so that $\tilde{f}$ extends $f$ as desired.

Finally, we prove that if $m>l-2$, then $\operatorname{im }d_{1}=\sHom(Q_{m} \otimes e_{m+2}A, e_{l}A)$.  To this end, given a right $A$-module morphism $f:Q_{m} \otimes e_{m+2}A \rightarrow e_{l}A$, we must show that it factors as $Q_{m} \otimes e_{m+2}A \rightarrow A_{m,m+1} \otimes e_{m+1}A \overset{g}{\rightarrow} e_{l}A$.  With the notation as in the previous paragraph, we note that $f(\sum_{p} x_{p} \otimes x_{p}^{*} \otimes 1)$ is in $A_{l,m+2}$ so has the form $\sum_{p}y_{p}x_{p}^{*}$ for some $y_{1}, \ldots, y_{n} \in A_{l, m+1}$.  We define $g$ by letting $g(x_{p} \otimes 1)=y_{p}$.
\end{proof}

\begin{cor} \label{cor.someext}
If $n \geq 1$ and $M$ is in ${\sf Gr }A$, then $\usExt^{i}(A/A_{\geq n},M)=0$ for $i>2$ and $\usExt^{i}(A/A_{\geq n},e_{l}A)=0$ for $i \neq 2$.
\end{cor}

\begin{proof}
We prove the result by induction on $n$.  When $n=1$, the first result follows from the remark immediately preceding Theorem \ref{thm.gorenstein}, while the second result follows from Theorem \ref{thm.gorenstein}.

For the general case, we note that, by Lemma \ref{lemma.basic2}(1), the exact sequence
$$
0 \rightarrow A_{\geq n}/A_{\geq n+1} \rightarrow A/A_{\geq n+1} \rightarrow A/A_{\geq n} \rightarrow 0
$$
in ${\sf Bimod }A$ induces a long exact sequence, of which
$$
\usExt^{i}(A/A_{\geq n}, M) \rightarrow \usExt^{i}(A/A_{\geq n+1},M) \rightarrow \usExt^{i}(A_{\geq n}/A_{\geq n+1},M)
$$
is a part.  If $i>2$ the left term is zero by induction while the right term is zero by Lemma \ref{lemma.extcommute} and induction.  Therefore, the center is zero in this case.  If $i=0$ or $i=1$ and $M=e_{l}A$, the same reasoning ensures that the center is zero.
\end{proof}

\begin{cor} \label{cor.gor0}
Suppose $M$ is an object of ${\sf Gr }A$.
\begin{enumerate}
\item{}
For $i>2$,
$$
\mbox{R}^{i}\tau (M)=0.
$$
\item{} For $i \neq 2$
$$
\mbox{R}^{i}\tau (e_{l}A)=0.
$$
\item{}
$$
(\mbox{R}^{2}\tau (e_{l}A))_{l-2-i} = \begin{cases} A^{*}_{l-2-i, l-2} & \text{if $i \geq 0$} \\ 0 & \text{otherwise} \end{cases}
$$

\end{enumerate}
\end{cor}

\begin{proof}
By Lemma \ref{lemma.torsfunct} and the fact that ${\sf Gr }A$ has exact direct limits, the first two results follow directly from Corollary \ref{cor.someext}.

To prove (3), we prove two preliminary results.  We first claim that
$$
\usExt^{2}(A/A_{\geq n+1}, e_{l}A) \subset [l-2-n,l-2].
$$
To this end, we note that, by Corollary \ref{cor.someext}, the sequence
$$
0 \rightarrow A_{\geq n}/A_{\geq n+1} \rightarrow A/A_{\geq n+1} \rightarrow A/ A_{\geq n} \rightarrow 0
$$
induces an exact sequence
\begin{equation} \label{eqn.ext2}
0 \rightarrow \usExt^{2}(A/ A_{\geq n}, e_{l}A) \rightarrow \usExt^{2}(A/A_{\geq n+1}, e_{l}A) \rightarrow \usExt^{2}(A_{\geq n}/A_{\geq n+1}, e_{l}A) \rightarrow 0.
\end{equation}
If $n=1$ in (\ref{eqn.ext2}), then the claim follows from Lemma \ref{lemma.extbound}, Theorem \ref{thm.gorenstein} and Lemma \ref{lemma.extcommute}.  The general case follows from the induction hypotheses and (\ref{eqn.ext2}).

We next claim that $\usExt^{2}(A/A_{\geq n+1}, e_{l}A)_{l-2-n} \cong A_{l-2-n, l-2}^{*}$.  To prove this, we note that when $n=0$, the claim follows from Theorem \ref{thm.gorenstein}.  For $n>0$,
\begin{eqnarray*}
\usExt^{2}(A/A_{\geq n+1}, e_{l}A)_{l-2-n} & \cong & \usExt^{2}(A_{\geq n}/A_{\geq n+1}, e_{l}A)_{l-2-n} \\
& \cong & A_{l-2-n, l-2}^{*}
\end{eqnarray*}
where the first isomorphism follows from the first claim and (\ref{eqn.ext2}), while the second isomorphism follows from the $n=0$ case and Lemma \ref{lemma.extcommute}.

Finally, we prove (3).  We have
\begin{eqnarray*}
(\mbox{R}^{2}\tau (e_{l}A))_{l-2-i} & \cong & \dlim \usExt^{2}(A/ A_{\geq n}, e_{l}A)_{l-2-i} \\
& \cong & \usExt^{2}(A/ A_{\geq i+1}, e_{l}A)_{l-2-i} \\
& \cong & A^{*}_{l-2-i, l-2},
\end{eqnarray*}
where the first isomorphism is from Lemma \ref{lemma.torsfunct}, the second isomorphism follows from the first claim, and the third isomorphism follows from the second claim.
\end{proof}

For the rest of Section \ref{section.gor}, in addition to our previous assumptions on $A$,
\begin{itemize}
\item{} we assume $A$ is coherent,

\item{} for each $i \in \mathbb{Z}$, we let $\mathcal{A}_{i}$ denote $\pi(e_{i}A) \in {\sf cohproj }A$, and

\item{} if $\underline{\iota}:{\sf cohproj }A \rightarrow {\sf Proj }A$ is the functor defined in Lemma \ref{lemma.faithful}, we abuse notation by letting $\mathcal{A}_{i}$ denote $\underline{\iota}(\pi(e_{i}A))$.
\end{itemize}

For the next two results, we use the fact that there are isomorphisms
\begin{eqnarray*}
\operatorname{Hom }_{{\sf Proj }A}(\mathcal{A}_{j}, \pi(-))  & \cong & \operatorname{Hom }_{{\sf Gr }A}(e_{j}A, \omega \pi(-)) \\
& \cong & (\omega \pi (-))_{j}.
\end{eqnarray*}
\begin{cor} \label{cor.gorenstein}
There are isomorphisms of $\operatorname{End }\mathcal{A}_{i}-\operatorname{End }\mathcal{A}_{j}$-bimodules
$$
\operatorname{Ext}_{{\sf cohproj }A}^{q}(\mathcal{A}_{j},\mathcal{A}_{i}) \cong \begin{cases} A_{ij} & \mbox{if $q=0$} \\ A^{*}_{j,i-2} & \mbox{if $q=1$} \end{cases}.
$$
\end{cor}

\begin{proof}
If $q=0,1$,  there are isomorphisms of $\operatorname{End }\mathcal{A}_{j}-\operatorname{End }\mathcal{A}_{i}$-bimodules
$$
\operatorname{Ext}_{{\sf cohproj }A}^{q}(\mathcal{A}_{j},\mathcal{A}_{i}) \cong \operatorname{Ext}_{{\sf Proj }A}^{q}(\mathcal{A}_{j},\mathcal{A}_{i})
$$
by Corollary \ref{cor.gor0} and Lemma \ref{lemma.faithful}.  In addition, by the remark preceding the corollary, for $q=0,1$ there are isomorphisms
$$
\operatorname{Ext}^{q}_{{\sf cohproj }A}(\mathcal{A}_{j}, \mathcal{A}_{i}) \cong (\mbox{R}^{q}\omega (\mathcal{A}_{i}))_{j}.
$$
Furthermore, by Corollary \ref{cor.gor0}(2), $\tau (e_{l}A)=\operatorname{R}^{1}\tau (e_{l}A)=0$.  Thus, by Theorem \ref{theorem.quotients},
$$
\operatorname{Hom}_{{\sf cohproj }A}(\mathcal{A}_{j}, \mathcal{A}_{i}) \cong A_{ij},
$$
and, by Theorem \ref{theorem.quotients} and Corollary \ref{cor.gor0}(3),
\begin{eqnarray*}
\operatorname{Ext}^{1}_{{\sf cohproj }A}(\mathcal{A}_{j}, \mathcal{A}_{i}) & \cong & (\operatorname{R}^{2}\tau (e_{i}A))_{j} \\
& \cong & A^{*}_{j, i-2}.
\end{eqnarray*}
\end{proof}

We now prove a generalization of Serre vanishing \cite[Theorem 3.5(2)]{finite}.

\begin{theorem} \label{thm.vanish}
If $M$ is a coherent object of ${\sf Gr }A$,
$$
\operatorname{Ext}^{1}_{{\sf Proj }A}(\mathcal{A}_{i}, \pi(M))=0
$$
whenever $i>>0$.
\end{theorem}

\begin{proof}
By the remark preceding Corollary \ref{cor.gorenstein}, there is an isomorphism
$$
\operatorname{Ext}^{1}_{{\sf Proj }A}(\mathcal{A}_{i}, \pi M) \rightarrow \mbox{R}^{1}\omega (\pi(M))_{i}.
$$
By Theorem \ref{theorem.quotients}, the right-hand side is isomorphic to $(\mbox{R}^{2}\tau M)_{i}$. Since $\mbox{R}^{3}\tau=0$ by Corollary \ref{cor.gor0}(1), it suffices, by considering a presentation of $M$, to show that, for any $j$, $(\mbox{R}^{2}\tau e_{j}A)_{i}=0$ for $i>>0$.  This follows immediately from Corollary \ref{cor.gor0}.
\end{proof}

\begin{lemma} \label{lemma.sofinite}
If $M$ is coherent, then $(\mbox{R}^{i}\tau M)_{j}$ and $(\mbox{R}^{i}\omega (\pi(M))_{j}$ are finite-dimensional over $A_{jj}$ for all $i \geq 0$.
\end{lemma}

\begin{proof}
We first show that $(\mbox{R}^{i}\tau M)_{j}$ is finite-dimensional over $A_{jj}$ for all $i \geq 0$.  By Corollary \ref{cor.gor0}(2) and (3), the result holds when $M \cong e_{i_{1}}A \oplus \cdots \oplus e_{i_{n}}A$.  Next, we note that $M$ coherent implies that there is a short exact sequence in ${\sf Gr }A$
$$
0 \rightarrow R \rightarrow \oplus_{l \in I}e_{l}A \rightarrow M \rightarrow 0
$$
where $I$ is finite and $R$ is coherent.  Therefore, by Corollary \ref{cor.gor0}(1), the first result follows from descending induction on $i$ as in the proof of \cite[Lemma 3.2]{finite}.

The second part of the lemma follows from the first part and Theorem \ref{theorem.quotients}.
\end{proof}
\subsection{Proof of Theorem \ref{theorem.supermain}(5)} \label{section.proof2}

\begin{theorem}
The sequence $(\mathcal{A}_{-i})_{i \in \mathbb{Z}}$ in ${\sf cohproj }A$ is a helix.
\end{theorem}

\begin{proof}
Properties (1), (2), (3), and (6) follow immediately from Corollary \ref{cor.gorenstein}, while property (7) follows from the fact that $A_{01}$ is 2-periodic.  Property (4) follows from (\ref{eqn.ei}) since the right-dimension of $A_{-i-2, -i-1}$ equals the left-dimension of $A_{-i-1, -i}$.  Property (5) follows from Proposition \ref{prop.phiiso}.

Next, we check property (8).   By the remark preceding Corollary \ref{cor.gorenstein}, there are $A_{jj}$-vector space isomorphisms
$$
\operatorname{Hom }_{{\sf Proj }A}(\mathcal{A}_{j},\pi M) \cong \operatorname{Hom }_{{\sf Gr }A}(e_{j}A,\omega \pi M) \cong (\omega \pi M)_{j}.
$$
Therefore, property (8) follows from Lemma \ref{lemma.faithful}(2) and Lemma \ref{lemma.sofinite}.

It remains to check property (9).  We first show the indicated sequence is projective in the sense of Section \ref{section.ampleness}.  To this end, we first claim that if $f:\mathcal{N} \rightarrow \mathcal{M}$ is an epimorphism in ${\sf cohproj }A$ and $\mathcal{N}$ is a finite direct sum of modules of the form $\mathcal{A}_{l}$, then there exists an $n$ such that $\operatorname{Hom}_{{\sf cohproj }A}(\mathcal{A}_{i},f)$ is surjective for all $i>n$.  By \cite[Corollaire 1, p. 368]{gab} and Corollary \ref{cor.gor0}(2), we may assume that $\operatorname{ker }f = \pi K$ where $\iota(K)$ is torsion-free.  By applying $\operatorname{Hom}_{{\sf cohproj }A}(\mathcal{A}_{i},-)$ to this sequence, the claim is reduced to showing that there exists an $n$ such that for all $i>n$, $\operatorname{Ext}^{1}_{{\sf cohproj }A}(\mathcal{A}_{i},\pi(K))=0$.  This last fact follows from Theorem \ref{thm.vanish} and Lemma \ref{lemma.faithful}(3).

For the general case, we use the fact that there is an epimorphism from a finite direct sum of modules of the form $\mathcal{A}_{i}$ to $\mathcal{N}$, which induces an epimorphism, $g$, to $\mathcal{M}$.  The claim then implies that there is some $n$ such that for all $i>n$, $\operatorname{Hom}_{{\sf cohproj }A}(\mathcal{A}_{i},g)$ is surjective, which then implies that $\operatorname{Hom}_{{\sf cohproj }A}(\mathcal{A}_{i},f)$ is surjective, as desired.

Finally, the fact that our sequence is ample follows from the definition of coherence and property (4).
\end{proof}

\section{An application} \label{section.piont}
In this section we confirm the $\mathbb{P}^{1}_{n}$ has a helix $\underline{\mathcal{L}}$.  Corollary \ref{cor.main} follows as an immediate consequence.  Throughout the section, we let
$$
A=k \langle x_{1}, \ldots, x_{n} \rangle/ (b),
$$
where $x_{i}$ has degree one for all $i$, and  $b=\sum_{i=1}^{n}x_{i}\sigma(x_{n-i})$ for some graded automorphism $\sigma$ of the free algebra.  Since $A$ is coherent by \cite[Theorem 4.3]{piont}, ${\sf cohproj }A$ is abelian, and ${\sf cohproj }A$ is $k$-linear by \cite[Proposition B8.1]{az2}.  For $j \in \mathbb{Z}$, we let $[j]$ denote the shift functor on ${\sf Gr }A$, so that $M[j]_{i} := M_{i+j}$.  We let $\pi:{\sf coh }A \rightarrow {\sf cohproj }A$ denote the quotient functor, and we define
$$
\mathcal{L}_{i} := \pi(A[i]).
$$
By \cite[Proposition 5.1(2)]{piont}, $\operatorname{Hom}(\mathcal{L}_{-1},\mathcal{L}_{0}) \cong A_{1}$, and is thus an $n$-dimensional vector space over $k$.  Furthermore, by the proof of \cite[Proposition 5.1(2)]{piont}, $\operatorname{End}\mathcal{L}_{i}=k$.  Therefore, by Corollary \ref{cor.central}, it suffices to confirm properties (2), (4), (6), (8) and (9) of $\underline{\mathcal{L}}$.  By the proof of \cite[Proposition 5.1(2)]{piont}, the first equality of (2) holds.  By \cite[Proposition 1.5(b)]{piont}, $\mathbb{P}^{1}_{n}$ satisfies Serre duality (see \cite[Section 1.3]{kgroup} for the exact form of the duality).  The second part of (2), as well as property (6), follows from duality.  Property (4) follows from \cite[Section 1.9]{kgroup}.  Property (8) follows from \cite[Proposition 5.1(3)]{piont}.  Finally, property (9) is observed at the end of the statement of \cite[Proposition 5.1]{piont}.

\end{document}